\documentclass[10pt]{amsart}
 \usepackage{amssymb,latexsym,amsmath,amsthm}

 \renewcommand{\div}{\mathop{\mathrm{div}}\nolimits}

 \usepackage[numbers,sort&compress]{natbib}
 
\usepackage{fullpage}

 \renewcommand{\div}{\mathop{\mathrm{div}}\nolimits}

\def\R{{\mathbb R} }

\def\r{\R^{n+1}_{+}}

\def\br{\partial\r}
\def\super{\overline}

\newtheorem{thm}{Theorem}[section]

\newtheorem{dfn}{Definition}[section]

\newtheorem{lemma}{Lemma}[section]
\newtheorem{notation}{Notation}[section]

\newtheorem{rem}{Remark}[section]

\numberwithin{equation}{section}

\begin{document}
\date{}

\title{Solutions of multi-component  fractional symmetric systems}

\author{Mostafa Fazly}

\address{Department of Mathematics, The University of Texas at San Antonio, San Antonio, TX 78249, USA}
\email{mostafa.fazly@utsa.edu}

\address{Department of Mathematics \& Computer Science, University of Lethbridge, Lethbridge, AB T1K 3M4 Canada.}
\email{mostafa.fazly@uleth.ca}

\thanks{The author gratefully acknowledges Natural Sciences and Engineering Research Council of Canada (NSERC) Discovery Grant and University of Texas at San Antonio Start-up Grant.}

\maketitle

\begin{abstract} We study the following elliptic system concerning the fractional Laplacian operator
$$(- \Delta)^ {s_i}  u_i = H_i ( u_1,\cdots,u_m) \ \ \text{in}\ \ \mathbb{R}^n,$$
when  $0<s_i<1$,  $u_i: \mathbb R^n\to R$ and $H_i$ belongs to 
$C^{1,\gamma}(\mathbb{R}^m)$ for $\gamma > \max(0,1-2\min \left \{s_i \right \})$ for $1\le i \le m$.  The above system is called symmetric when the matrix $\mathcal H=(\partial_j H_i(u_1,\cdots,u_m))_{i,j=1}^m$ is symmetric. The notion of symmetric systems seems crucial to study this system with a general nonlinearity $H=(H_i)_{i=1}^m$.  We establish De Giorgi type results for stable and $H$-monotone solutions of symmetric systems in lower dimensions that is either  $n=2$ and $0<s_i<1$  or $n=3$ and $1/2 \le \min\{s_i\}<1$. The case that $n=3$ and at least one of parameters $s_i$ belongs to $(0,1/2)$ remains open as well as the case  $n \ge 4$.   
 Applying a geometric Poincar\'{e} inequality, we  conclude that  gradients of components of solutions are parallel in lower dimensions when the  system is coupled. More precisely, we show that the angle between vectors $\nabla u_i$ and $\nabla u_j$ is exactly $\arccos\left({|\partial_j H_i(u)|}/{\partial_j H_i(u)}\right)$.   
 In addition, we provide  Hamiltonian identities, monotonicity formulae and  Liouville theorems.  Lastly, we  apply some of our main results to    a two-component nonlinear Schr\"{o}dinger system, that is a particular case of the above system, and we prove Liouville theorems and monotonicity formulae.

\end{abstract}

\noindent
{\it \footnotesize 2010 Mathematics Subject Classification}. {\scriptsize  35J60, 35J50, 35B35, 35B45}\\
{\it \footnotesize Keywords:  Nonlinear elliptic systems,  fractional Laplacian operator,  Hamiltonian identity, monotonicity formula, symmetry of entire solutions}. {\scriptsize }

%\tableofcontents

\section{Introduction}
We examine the following  system of nonlocal elliptic equations for $1\le i \le m$, 
\begin{eqnarray}
\label{mainh}
(- \Delta)^ {s_i}  u_i = H_i ( u_1,\cdots,u_m) \ \ \text{in}\ \ \mathbb{R}^n,
  \end{eqnarray}
when  $(-\Delta)^{s_i}$ stands for the fractional Laplacian operator and $ u_i:\mathbb{R}^n\to \mathbb{R}$. We assume that each parameter $s_i$ belongs to $(0,1)$ and  $H=(H_i)_{i=1}^m$ is a sequence of functions that each component  $H_i$ belongs to $C^{1,\gamma}(\mathbb{R}^m)$ for $\gamma > \max(0,1-2\min \left \{s_i \right \})$.  We also assume that the nonlinearity $H $ satisfies 
 \begin{equation}
\partial_i H_j (u)\partial_j H_i (u)> 0 \ \ \ \text{when} \ \ \ \partial_j H_i(u):=\frac{\partial H_i(u)}{\partial{u_j}} \ \ \   \text{for} \ \ \ 1\le i<j\le m. 
 \end{equation}
  More recently,  the fractional Laplacian operator $(-\Delta)^{s_i}$ has been of great interests in the literature and various properties of the operator are explored.  There are different mathematical approaches to  define this operator.  The fractional Laplacian operator on $\mathbb R^n$ can be defined as a pseudo-differential operator using the Fourier transform 
 \begin{equation}
\widehat{ (-\Delta)^{s_i} w}(\zeta)=|\zeta|^{2s_i} \widehat w(\zeta), 
 \end{equation}
 when the hat operator is given  by 
 \begin{equation}
 \widehat w(\zeta)=(2\pi)^{-n/2} \int_{\mathbb R^n} w(x) e^{-i\zeta\cdot x} dx.
 \end{equation}
 Suppose that each  $u_i\in C^{2\sigma}(\mathbb R^n)$ for $\sigma >s_i>0$ and 
$ \int_{\mathbb R^n} \frac{|u_i(z)|}{(1+|z|)^{n+2s_i}} dz<\infty.$ Then,  the fractional Laplacian of $u_i$ can be also defined as
\begin{equation}
(-\Delta)^{s_i} u_i(x):= P.V. \int_{\mathbb R^n} \frac{u_i(x)-u_i(z)}{|x-z|^{n+2s_i}} dz ,
\end{equation}
 for every $x\in\mathbb R^n$ where P.V. stands for the principal value.   It is by now standard that the fractional Laplacian $(-\Delta )^s$ operator  where $s$ is any positive, noninteger number can be denoted  as the
Dirichlet-to-Neumann map for an extension function satisfying a higher order elliptic equation in the upper half space with one extra spatial dimension, see \cite{cafs} by Caffarelli and Silvestre.  To be mathematically more precise,  let $s_i\in(0,1)$ and $u_i\in C^{2\sigma}(\mathbb R^n)\cap L^1\left(\mathbb R^n,(1+\vert z\vert)^{n+2s_i}dz\right)$ when $\sigma>s_i$ for all  $i=1,\cdots,m$. For $\bar x=(x,y)$ in $\R^{n+1}_{+}$,  set
 \begin{equation}
 v_i(\bar x)= \int_{\mathbb R^n} P_i(x-z ,y )u_i(z)\;dz \ \ \ \text{when} \ \ \  
P_i(x,y) = p_{n,s_i} \,\frac{y^{2s_i}}{  \left( |x|^2+|y|^2 \right)^{\frac{n+2s_i}{2}} } , 
 \end{equation}
and $p_{n,s_i}$ is a normalizing constant.   Then, $v=(v_i)_{i=1}^m$ satisfies the following system 
\begin{equation}
\left\{
\begin{aligned}
\div(y^{1-2s_i}\nabla  v_i)&=0&\quad\text{in \quad $\mathbb R^{n+1}_{+}$}\\
 v_i&= u_i&\quad\text{on $\br$,}\\
-\lim_{y\to0} y^{1-2s_i}\partial_{y} v_i&= d_{s_i} (-\Delta)^{s_i} u_i&\quad\text{on $\br$,}
\end{aligned}
\right.
\end{equation}
for the  constant 
\begin{equation} \label{kappas}
d_{s_i} = \frac{\Gamma(1-s_i)}{2^{2s_i-1} \Gamma(s_i)},
\end{equation} 
and  $ v_i\in C^2(\r)\cap C(\super\r)$ and $y^{1-2s_i}\partial_{y} v_i\in C(\super\r)$.   In short, the extension function $v=(v_i)_{i=1}^m$  satisfies the system 
 \begin{eqnarray}\label{emainh}
 \left\{ \begin{array}{lcl}
\hfill \div(y^{a_i} \nabla v_i)&=& 0   \ \ \text{in}\ \ \mathbb{R}_+^{n+1}=\left \{x \in \mathbb R^n, y>0 \right \},\\   
\hfill -\lim_{y\to0}y^{a_i} \partial_{y} v_i&=& d_{s_i} H_i( v)   \ \ \text{in}\ \ \partial\mathbb{R}_+^{n+1}=\left \{x \in \mathbb R^n, y=0 \right \},
\end{array}\right.
  \end{eqnarray}
when $a_i=1-2s_i$.    For the sake of simplicity,  we fix the  following notation. %before stating various definitions.  
\begin{notation}
Suppose that $0<s_i<1$ for $1\le i\le m$. Then, $s_*=\min_{1\le i\le m}\{s_i\}$ and  $s^*=\max_{1\le i\le m}\{s_i\}$.   
\end{notation}

We now define the notion of stable solutions by linearizing system (\ref{emainh}). 
 
\begin{dfn}  A solution $ v=(v_i)_{i=1}^m$ of  system (\ref{emainh}) in $\mathbb{R}_+^{n+1}$ is called  stable, if there exists a sequence $\phi=(\phi_i)_{i=1}^m$ in $C^\infty(\mathbb{R}_+^{n+1})$ satisfying
\begin{eqnarray}\label{stabilityh}
 \left\{ \begin{array}{lcl}
\hfill \div(y^{a_i} \nabla \phi_i)&=& 0   \ \ \text{in}\ \ \mathbb{R}_+^{n+1},\\   
\hfill -\lim_{y\to0}y^{a_i} \partial_{y} \phi_i&=& d_{s_i} \sum_{j=1}^m \partial_{j}H_i ( v)\phi_j   \ \ \text{in}\ \ \partial\mathbb{R}_+^{n+1},
\end{array}\right.
  \end{eqnarray}
 when $\partial_{j}H_i ( v)\phi_j \phi_i>0$ for all $i,j$  and each $\phi_i $ does not change sign  for all $i=1,\dots,m$. 
 \end{dfn}
 
  In 1978, De Giorgi \cite{DeGiorgi} conjectured that bounded monotone solutions of the Allen-Cahn equation 
 \begin{equation}\label{AllenCahn}
 \Delta u+u-u^3=0 \ \ \ \text{in} \ \ \ \mathbb{R}^n,
 \end{equation}
  are one-dimensional solutions at least up to eight dimensions. This conjecture is known to be true for dimensions $n=2,3$ by Ghoussoub-Gui in \cite{gg} and Ambrosio-Cabr\'{e} in \cite{ac}, respectively.  There is an example by del Pino, Kowalczyk and Wei in \cite{dkw} that shows the dimension eight is the critical dimension. For dimensions $4\le n \le 8$ there are various partial results under certain extra (natural) assumptions on solutions by Ghoussoub and Gui in \cite{gg2},  Savin in \cite{sav} and references therein. The remarkable point is that in lower dimensions,  that is when $n\le 3$,  this conjecture holds for  the scalar equation 
   \begin{equation}\label{sceqH}
-\Delta u=H(u) \ \ \ \text{in} \ \ \ \mathbb{R}^n,
 \end{equation}
   with a general nonlinearity $H\in C^1(\mathbb R)$. In this regard, we refer interested readers to \cite{gg} by   Ghoussoub and Gui for $n=2$ and to \cite{ac} by Alberti, Ambrosio and Cabr\'{e} for $n=3$ and also to  \cite{fsv} by Farina, Sciunzi and Valdinoci  for a geometric approach in two dimensions.    
  When the Laplacian operator in (\ref{sceqH}) is replaced with the fractional Laplacian operator $ (-\Delta)^s$,  De Giorgi type results are known  when either $n=2$ and $0<s<1$ or $n=3$ and $1/2\le s<1$, see \cite{sv,cs1,cs2,cso}.  However, in three dimensions and when  $0<s<1/2$  the problem seems to be more challenging and remains open.   
  
More recently,  Ghoussoub and the author in \cite{fg} provided De Giorgi type results up to three dimensions  for the following  elliptic gradient systems  when $u:\mathbb R^n\to \mathbb R^m$ 
    \begin{equation}\label{syeqH}
  -\Delta u=\nabla H(u)
  \ \ \ \text{in} \ \ \ \mathbb{R}^n,
 \end{equation}
  for a general nonlinearity $H\in C^2(\mathbb R^m)$. In this regard,  authors introduced the notions  of orientable systems and $H$-monotone solutions to adjust and to apply the mathematical techniques and ideas given for the scalar equation (\ref{sceqH}) to (\ref{syeqH}).  It seems that these concepts are essential, in this context,  to explore system of equations.  Note also that Sire and the author in \cite{fs} studied  system (\ref{syeqH}),  when the Laplacian operator is replaced with the fractional Laplacian operator $ (-\Delta)^s$,  and provided De Giorgi type results  for certain parameters $s=(s_1,\cdots,s_n)\in\mathbb R^n$ in lower dimensions.     
   \begin{dfn} \label{weak} 
System  (\ref{emainh}) is called orientable, if there exist nonzero functions $\theta_k\in C^1(\mathbb{R}^{n+1}_+)$, $k=1,\cdots,m$, which do not change sign, such that for all $i,j$ with $1\leq  i<j\leq m$,
 \begin{equation}\label{oriantableu}
 \hbox{$ \partial_{j}H_i({ u}) \theta_i( x)\theta_j( x) > 0$ \, for all $ x\in\mathbb{R}^n$.}
  \end{equation} 
  Similarly, if the condition (\ref{oriantableu}) holds for the extension $ v$, then we say (\ref{emainh}) is orientable. 
\end{dfn}
Note that the orientability imposes a combinatorial assumption on the sign of the nonlinearity $H=(H_i)_{i=1}^m$ and therefore on the system.  

\begin{dfn} 
We say that a solution ${u}=(u_i)_{i=1}^m$ of (\ref{mainh}) is $H$-monotone if the following hold,
\begin{enumerate}
 \item For every $i\in \{1,\cdots, m\}$, $u_i$ is strictly monotone in the $x_n$-variable (i.e., $\partial_n u_i\neq 0$).

\item  For $i<j$, 
we have 
  \begin{equation}
\hbox{$\partial_{j}H_i({ u}) \partial_n u_i( x) \partial_n u_j ( x) > 0$  for all $x\in\mathbb{R}^n$.}
\end{equation}

\end{enumerate}
Similarly, if above conditions  hold for the extension function $ v$ then we say that $ v$ is $H$-monotone.  
\end{dfn}
    
It is straightforward to observe that $H$-monotonicity implies stability. One can show this by differentiating (\ref{mainh}) with respect to the variable $x_n$ and setting $\phi_i:=\partial_n u_i$.   As the last definition of this section,  we provide the notion of  symmetric systems that plays  a key role in our main results. 

\begin{dfn}\label{symmetric} We call system (\ref{mainh}) symmetric if the matrix of gradient of all components of $H$ that is 
\begin{equation}
\mathbb{H}:=(\partial_j H_i(u))_{i,j=1}^{m}
\end{equation}
  is symmetric.  Similarly, if $\mathbb{H}:=(\partial_j H_i(v))_{i,j=1}^{m}$ is symmetric for the extension $ v$, then we say (\ref{emainh}) is symmetric. 
 \end{dfn}

 In this article,  we prove De Giorgi type results for stable and $H$-monotone solutions of symmetric system (\ref{emainh}) with a general nonlinearity $H$ when  $n=2$ and $0<s_i<1$  or $n=3$ and $1/2 \le s_*<1$, see Theorem \ref{thsymv}.  We also provide Liouville theorems for stable solutions of  (\ref{emainh})  under some extra assumptions on the nonlinearity $H$ when $n\le 2(1+s_*)$, see Theorem \ref{lioupositive}.  In addition, we establish the following geometric Poincar\'{e} inequality for bounded stable solutions of (\ref{emainh}). Let $\eta=(\eta_k)_{k=1}^m \in C_c^1(\mathbb R_+^{n+1})$, then 
\begin{eqnarray}
&&\sum_{i=1}^m  \frac{1}{d_{s_i}}   \int_{    \{|\nabla_{ x}  v_i|\neq 0 \} \cap \mathbb R_+^{n+1} }      y^{1-2s_i} \left(   |\nabla v_i|^2 \mathcal{A}_i^2 + | \nabla_{T_i} |\nabla_{ x}  v_i| |^2  \right)\eta_i^2 d\bar x
\\&& \label{secondH} +\sum_{i\neq j} \int_{\partial\mathbb R_+^{n+1}}   \left(  \sqrt{\partial_{j} H_i( v) \partial_{ i} H_j( v) }  |\nabla_{ x}  v_i|  |\nabla_{ x}  v_j| \eta_i \eta_j   - \partial_{j}H_i(v)  \nabla_{ x}  v_i \cdot   \nabla_{ x}  v_j \eta_i^2 \right) dx
\\&\le& \sum_{i=1}^m \frac{1}{d_{s_i}}  \int_{\mathbb R_+^{n+1}} y^{1-2s_i}  |\nabla_{ x}  v_i|^2   |\nabla \eta_i|^2 d\bar x,
  \end{eqnarray} 
 where $\nabla_{T_i}$ stands for the tangential gradient along a given level set of $v_i$ and 
$\mathcal{A}_i^2$ for the sum of the squares of the principal curvatures of such a level set, see Theorem \ref{lempoin}. Note that for the case of scalar equations, $m=1$,  a similar inequality was established by Sternberg and  Zumbrun \cite{sz} to study phase transitions and area-minimizing surfaces.  For this case,  the boundary term (\ref{secondH}) disappears.   The idea of applying the  geometric Poincar\'{e} inequality  to prove De Giorgi type results  was initiated by Farina,  Sciunzi,  Valdinoci in \cite{fsv} and references therein. Ghoussoub and the author in \cite{fg} established a counterpart of this inequality for local systems of the form  (\ref{syeqH}).  Note that for  symmetric systems, $m\ge 2$, the boundary term (\ref{secondH})  becomes 
\begin{equation}
 \sum_{i\neq j} \int_{\partial\mathbb R_+^{n+1}}  \partial_{j} H_i( v)  \left(  |\nabla_{ x}  v_i|  |\nabla_{ x}  v_j| \eta_i \eta_j   -  \nabla_{ x}  v_i \cdot   \nabla_{ x}  v_j \eta_i^2 \right) dx,
 \end{equation}
where the integrand has a fixed sign for an appropriate test function $\eta=(\eta_i)_{i=1}^m$.  In the light of this inequality,  we prove De Giorgi type results in two dimensions and we show that vectors $\nabla_x v_i(x,0)$ and $\nabla_x v_j (x,0)$ for $i\neq j$ are parallel and the angle between two vectors is $\arccos\left(\frac{|\partial_j H_i(v)|}{\partial_j H_i(v)}\right)$.

   Consider the scalar equation (\ref{sceqH}) when $\tilde H \ge 0$ for $\tilde H'(u)=-H(u)$. Modica in \cite{mod} proved that the following pointwise estimate holds for bounded solutions 
   \begin{equation}
|\nabla u|^2 \le 2 \tilde H(u) \ \ \ \text{in} \ \ \ \mathbb R^n.
   \end{equation}
This inequality has been used in the literature to study entire solutions of semilinear elliptic equations and,  in particular,  to establish De Giorgi type results, see  \cite{aac,ac,bcn,dkw,gg,gg2,gui,sav}. Unfortunately, a counterpart of this inequality does not hold for the system of equations of the form (\ref{syeqH}) with a general nonlinearity.   However,  assuming that $n=1$ and $m\ge 1$ and multiplying the $i^{th}$ equation of (\ref{syeqH}) with  $u_i'$ and integrating,  it is straightforward to observe that the following Hamiltonian identity holds
\begin{equation}
\frac{1}{2}\sum_{i=1}^m |u_i'(x)|^2 +H(u(x)) \equiv C \ \ \text{for}\ \ x\in\mathbb R , 
\end{equation}
  when $C$ is a constant.   Gui in \cite{gui} %asked and answered whether there is a counterpart of this identity for partial differential equations and systems.   He 
examined system (\ref{syeqH}) when $n\ge1$ and $m\ge 1$ and proved the following elegant Hamiltonian identity in higher dimensions
\begin{equation}\label{ghamH}
 \int_{R^{n-1}} \left[ \sum_{i=1}^m \frac{1}{2} \left( |\nabla_{x'} u_i(x)|^2 - |\partial_n u_i(x)|^2 \right) - H(u(x)) \right] dx'\equiv C \ \ \text{for}\ \ x_n\in\mathbb R ,
 \end{equation}
 where $x=(x',x_n)\in\mathbb R^n$ and $C$ is a constant.    In this paper,  as Theorem \ref{hamiltonthm}, we provide a counterpart of (\ref{ghamH}) for the fractional system (\ref{emainh}) in one dimension that is %$n=1$ of the form 
 \begin{equation}
 \sum_{i=1}^{m} \frac{1}{2 d_{ s_i}} \int_0^\infty y^{1-2s_i} \left[  (\partial_x v_i)^2 - (\partial_y v_i)^2 \right] dy -  \tilde H(v(x,0)) \equiv C \ \ \text{for} \ \ \ x \in \mathbb R,
\end{equation}
when $C$ is a constant and $\partial_i \tilde H(v)= H_i(v)$ for every $i$. Proving a similar identity in higher dimensions, $n\ge 2$, remains an open problem.  Our methods and ideas are strongly motived by the ones provided by  Gui in \cite{gui}, Cabr\'{e} and Sol\'{a}-Morales in \cite{cso}, Cabr\'{e} and Sire in \cite{cs1,cs2} and Sire and the author in \cite{fs}.

Lastly, we consider a  two-component system of the form (\ref{mainh}) with the following particular nonlinearity 
\begin{equation}
H_1(u_1,u_2)=\mu_1 u_1^3+\beta u_2^2 u_1 \ \ \text{and} \ \  H_2(u_1,u_2)=\mu_2 u_2^3+\beta u_1^2 u_2.
\end{equation} 
where $u:\mathbb R^n\to\mathbb R^2$ and $\mu_1,\mu_2,\beta$ are constants.  Both local and nonlocal systems of this type, known as the nonlinear Schr\"{o}dinger system, have been studied extensively in the literature. For more information, we refer interested readers to \cite{wwe2,wwe,ww,tvz,nttv,lw,bdw,fn,dww} and references therein. The nonlinear Schr\"{o}dinger system is a natural counterpart of the following nonlinear Schr\"{o}dinger equation, 
\begin{equation}\label{NLSch}
 \Delta u-u+u^3=0\ \ \text{in}  \ \ \mathbb R^n. 
 \end{equation} 
Even though equations (\ref{NLSch}) and (\ref{AllenCahn})  look alike,  %with just a sign change in (\ref{AllenCahn}) this equation gives. Even though these equations look alike 
their solutions behave very differently.  In this article, we provide various Liouville theorems and monotonicity formulas for solutions of the nonlinear fractional Schr\"{o}dinger system, under some assumptions on parameters $s_1,s_2,\mu_1,\mu_2,\beta$.  

%Note that what we have mainly discussed so far refers to (\ref{main}) with a general nonlinearity, and the motivation for this was the fact that De Giorgi's conjecture for the Allen-Cahn equation, that is $\Delta u+u-u^3=0$,  holds for a general nonlinearity in lower dimensions. Just a sign change in this equation gives $\Delta u-u+u^3=0$ that is known as the nonlinear Schr\"{o}dinger equation.  Even though these equations look alike their solutions  have very different behaviours.  In this paper we shall study the Schr\"{o}dinger system with two components that is 
%\begin{eqnarray}\label{} \left\{ \begin{array}{lcl} \hfill (- \Delta)^s  u &=& \mu_1 u^3+\beta v^2 u   \ \ \text{in}\ \ \mathbb{R}^n,\\    \hfill (- \Delta)^t v &=& \mu_2 v^3+\beta u^2 v     \ \ \text{in}\ \ \mathbb{R}^n, \end{array}\right. \end{eqnarray}
%where $0<s,t<1$ and $\mu_1,\mu_2,\beta$ are constants. Note that this is a particular case of (\ref{mainh}). We shall provide Liouville theorems and monotonicity formulas for solutions of this system for various parameters $s,t,\mu_1,\mu_2,\beta$.  

%Finally let us mention that in the absence of monotonicity and stability assumptions the qualitative behaviour of solutions of elliptic and Hamiltonian systems are studies extensively in the literature. We shall refer interested readers to \cite{dff} by De Figueiredo and Felmer and to \cite{dfd} by De Figueiredo and Ding and references therein. 

The organization of the paper is as follows.   In the next section, we  provide some standard regularity results and estimates regarding fractional Laplacian operator. In Section \ref{decay}, we prove various technical tools needed to establish our main results. We provide  a Hamiltonian identity,  a geometric Poincar\'{e} inequality  and various gradient and energy estimates for system (\ref{emainh}) with a general nonlinearity $H=(H_i)_{i=1}^m$. We also provide  a monotonicity formula to conclude the optimality of  gradient and energy estimates.      In Section \ref{nonlin},  we establish De Giorgi type results and Liouville theorems for symmetric systems in lower dimensions. In addition, applying the Hamiltonian identity we analyze directional derivatives of the extension function $v$ satisfying (\ref{emainh}).   Lastly, in Section \ref{sch}, we consider a two-component fractional Schr\"{o}dinger system,  that is a particular case of (\ref{mainh}), and we provide various Liouville theorems and monotonicity formulae.     
    
  %   some decay estimates on solutions and a geometirc Poincar\'{e} inequality needed in our proofs of the main results.  In Section \ref{nonlin} we provide De Giorgi type results and Liouville theorems for symmetric systems as well as a Hamiltonian identity for system (\ref{emainh}) with a general nonlinearity $H=(H_i)_{i=1}^m$.    Finally, in Section \ref{sch} we provide Liouville theorems and monotonicity formulas for a particular system that is called the fractional Schr\"{o}dinger system. 

 %%%%%%%%%%%%%%%%%%%%%%%%%%%%%%%%%%%%%%%%
\section{Standard Elliptic Estimates}\label{pre}
%%%%%%%%%%%%%%%%%%%%%%%%%%%%%%%%%%%%%%
In this section,  we provide some standard estimates regarding fractional Laplacian operator.   We assume that  the nonlinearity $H=(H_i)_{i=1}^m$ and each $H_i$ belongs to $C^{1,\gamma}(\mathbb R^m)$  with $\gamma >\max(0,1-2s_*)$ and $ u=(u_i)_{i=1}^m$ is a sequence of bounded functions. We omit the proofs of   lemmata  in this section and we refer interested readers to \cite{cs1,cs2,ccinti,ccinti1,fs,sv} and references therein.  We start with the following regularity result for solutions of (\ref{mainh}).

\begin{lemma}\label{regu}
Suppose that $u=(u_i)_i$ is a bounded solution of  \eqref{mainh}. Then,  each $u_i$ is $C^{2,\beta}(\mathbb R^n)$ for some $0 <\beta < 1$ depending only on $s_i$ and $\gamma$.
\end{lemma}

For the sake of convenience,  we use the following notation used frequently in the literature.%,% see \cite{ccinti,cs1,cs2,cafs,sv} and references therein.   
\begin{notation}
Set $B_R^+=\{\bar x=(x,y)\in\mathbb R_+^{n+1}, |\bar x|<R\}$, $\partial^+ B_R^+=\partial B^+_R\cap \{y>0\}$,  $\Gamma^0_R=\partial B^+_R\cap \{y=0\}$ and  $C_R=B_R\times (0,R)$ for $R>0$.  
\end{notation}
We now provide a  regularity result for solutions of an equation in the half-space.  This can be applied to extensions functions $v$ satisfying (\ref{emainh}). 

\begin{lemma}\label{regv}
 Let $f_i \in C^\sigma (\Gamma^0_{2R})$ for some $\sigma \in (0,1)$, $R>0$ and 
$v_i \in L^\infty(B^+_{2R}) \cap H^1(B^+_{2R},y^{1-2s_i})$ for some $s_i \in (0,1)$ be a weak solution of
%\begin{equation*} \label{problemBR} \begin{cases} \div(y^{1-2s_i}\nabla v_i)=0&\text{ in } B^+_{2R}\subset\mathbb R^{n+1}_+\\  -\lim_{y \to 0} y^{1-2s_i}\frac{\partial v_i}{\partial y} =f_i  &\text{ on } \Gamma^0_{2R}. \end{cases} \end{equation*}
\begin{eqnarray}\label{}
 \left\{ \begin{array}{lcl}
\hfill \div(y^{1-2s_i}\nabla v_i)&=& 0   \ \ \text{in}\ \ B^+_{2R}\subset\mathbb R^{n+1}_+,\\   
\hfill -\lim_{y \to 0} y^{1-2s_i}\frac{\partial v_i}{\partial y} &=& f_i   \ \ \text{on}\ \ \Gamma^0_{2R}. 
\end{array}\right.
  \end{eqnarray}
Then, there exists  $\beta \in (0,1)$ depending only on $n$, $s_i$, and $\sigma$, 
such that $v_i \in C^{0,\beta}(\overline{B_R^+})$ and 
$y^{1-2s_i}\partial_y v_i \in C^{0,\beta}(\overline{B_R^+})$.  Furthermore, there exist constants $C_R$ and $D_R$ such that 
\begin{equation}
\|v_i\|_{C^{0,\beta}(\overline{B_R^+})} \leq C_i(R)\ \ 
\text{and} \ \ \|y^{1-2s_i} \partial_y v_i\|_{C^{0,\beta}(\overline{B_R^+})} \leq D_i(R),
\end{equation}
where $C_i(R)$ only depends on $n,a_i,R,|| v_i||_{L^\infty}(B^+_{2R})$ and $||f_i||_{L^\infty(\Gamma^0_{2R})}$ and $D_i(R)$ only depends on $n$, $a_i$, $R$, $|| v_i||_{L^\infty}(B^+_{2R})$ and $||f_i||_{C^\sigma(\Gamma^0_{2R})}$.
\end{lemma}

 The next lemma provides  gradient estimates for bounded solutions of \eqref{emainh} and it is a consequence of Lemma \ref{regu} and lemma \ref{regv}. For more information, we refer interested readers to \cite{LS, ccinti,cs1,cs2} and references therein.  We apply these estimates frequently in the proofs of our main results. 
  
\begin{lemma}\label{asymp}
Let $v=(v_i)_{i=1}^m$ be a bounded solution of \eqref{emainh}. Then,  each $v_i$ satisfies 
 \begin{eqnarray}
|\nabla_{ x} v_i( x,y)| &\le& C \ \ \text{for} \ \   x\in\mathbb{R}^{n} \ \ \text{and}\ \  y\ge0 , \\
|\nabla v_i ( x,y) | &\le& \frac{C}{1+y} \ \ \text{for} \ \ ( x,y)\in\mathbb {R}^{n+1}_+ , \\
%|\partial_y  v_i| &\le& \frac{C}{y}  \ \ \text{for} \ \  x\in\mathbb{R}^{n} \ \ \text{and}\ \ y>0 , \\
| y^{1-2s_i} \partial_y v_i ( x,y) | &\le& C \ \ \text{for} \ \ ( x,y)\in\mathbb {R}^{n+1}_+  , 
  \end{eqnarray}
  where the positive constant $C$ is independent from $x$ and $y$. 
\end{lemma}

 \section{Hamiltonian Identity and  Analytic and Geometric Estimates}\label{decay} 
 
We start this section with proving a Hamiltonian identity for solutions of (\ref{emainh}) in one dimension. We then apply this identity, see Theorem \ref{lioupositive},  to study  bounded stable solutions of symmetric (\ref{emainh}) with a general nonlinearity $H=(H_i)_{i=1}^m$.   Hamiltonian identities are amongst the most important tools to study qualitative behaviour of entire solutions of differential equations and in some cases they lead to properties such as monotonicity formulae.  

\begin{thm}\label{hamiltonthm}
Suppose that $v=(v_i)_{i=1}^m$ is a solution of (\ref{emainh}) in one dimension and $0<s_i<1$ for $1\le i\le m$.  Then the following Hamiltonian identity holds 
\begin{equation}\label{hamiltonianiden}
 \sum_{i=1}^{m} \frac{1}{2 d_{ s_i}} \int_0^\infty y^{1-2s_i} \left[  (\partial_x v_i)^2 - (\partial_y v_i)^2 \right] dy -  \tilde H(v(x,0)) \equiv C  \ \ \text{for} \ \ x \in \mathbb R,
\end{equation}
where $C$ is a constant that is independent from $x$ and $\partial_i \tilde H(v)= H_i(v)$ for every $i$.  
\end{thm}
%\begin{remark} Note that the concept of symmetric systems is crucial for the existence of such $\tilde H$.   \end{remark}
 \noindent\textbf{Proof:} Let $v=(v_i)_{i=1}^m$ be a solution of the extension problem (\ref{emainh}). Set 
\begin{equation}
w(x):=  \sum_{i=1}^{m}  \frac{1}{2 d_{s_i}}  \int_0^\infty y^{1-2s_i} \left[  (\partial_x v_i)^2 - (\partial_y v_i)^2 \right] dy.
\end{equation}
Differentiating  $w$ in terms of $x$, we get  
\begin{equation}\label{xw}
\partial_x w(x):= \sum_{i=1}^{m}  \frac{1}{d_{s_i}}  \int_0^\infty y^{1-2s_i} \left[  \partial_x v_i \partial_{xx} v_i - \partial_y v_i \partial_{xy} v_i \right] dy.
\end{equation}
From (\ref{emainh}) for each $i=1,\cdots, m$,  we have 
\begin{equation}
 y^{1-2s_i} \partial_{xx} v_i + \partial_y \left(   y^{1-2s_i}   \partial_y v_i  \right)=0. 
 \end{equation}
Combining this and (\ref{xw}), we end up with
\begin{equation}\label{xwi}
\partial_x w(x):= \sum_{i=1}^{m}  \frac{1}{d_{s_i}}   \int_0^\infty \left[  -\partial_x v_i \partial_y \left(   y^{1-2s_i}  \partial_y v_i  \right) -  y^{1-2s_i}  \partial_y v_i \partial_{xy} v_i \right] dy.
\end{equation}
Integration by parts for the first term in the above yields 
\begin{equation}
-\int_0^\infty \partial_x v_i \partial_y \left(    y^{1-2s_i}   \partial_y v_i  \right) dy = \int_0^\infty  y^{1-2s_i}   \partial_{xy} v_i  \partial_y v_i  dy  -  \lim_{y\to 0}  y^{1-2s_i}   \partial_{y} v_i    \partial_{x} v_i . 
\end{equation}
Substituting this in (\ref{xwi}) we get
 \begin{equation}
\partial_x w(x) = - \sum_{i=1}^{m}  \frac{1}{d_{s_i}}  \lim_{y\to 0}  y^{1-2s_i}    \partial_{x} v_i \partial_{y} v_i .
\end{equation}
From this and the boundary term in (\ref{emainh}),  we have 
 \begin{equation}
\partial_x w(x) = \sum_{i=1}^{m}  \partial_{i}  H(v(x,0)) \partial_{x} v_i = \partial_{x} \left(  \tilde H(v(x,0))  \right) .
\end{equation}
Therefore, 
 \begin{equation}
 \partial_x\left[    w(x) - \tilde H(v(x,0))   \right]=0.\end{equation}
This implies that $w(x) - \tilde H(v(x,0)) $ is constant in terms of $x$.

               \hfill $ \Box$

We now prove an inequality, known as the stability inequality, for stable solutions of (\ref{emainh}). This inequality plays an important role in this paper. 

%\subsection{Geometrical Inequality for stable solutions}
%In this section, we provide a few technical elliptic estimates needed in our proofs. Let us start with the following stability inequality that plays an important role in this paper. 

\begin{lemma}\label{stabilitylem} 
Let $ v=(v_i)_{i=1}^m$ be a stable solution of system (\ref{emainh}). Then,  the following stability inequality holds 
\begin{equation}\label{stabilityinh}
\sum_{i,j=1}^{m} \int_{\partial\mathbb{R}^{n+1}_+}  \sqrt{d_{s_i}d_{s_j}  \partial_{j}H_i( v) \partial_{v_i}H_j(v) }   \zeta_i \zeta_j  dx  \le \sum_{i=1}^{m} \int_{\mathbb {R}^{n+1}_+}  y^{a_i} |\nabla \zeta_i|^2  d\bar x , 
\end{equation}
for all $\zeta_i\in C_c^1( \mathbb R_+^{n+1})$.  
\end{lemma}

\noindent  \textbf{Proof:} Since  $ v=(v_i)_{i=1}^m$ is a stable solution of system (\ref{emainh}), there exists a sequence $\mathbf \phi=(\phi_i)_{i=1}^m$ that satisfies (\ref{stabilityh}).  Multiply the $i^{th}$ equation of (\ref{stabilityh}) with $\frac{\zeta_i^2}{\phi_i}$ and do integration by parts to get 
\begin{equation} - \int_{\mathbb R_+^{n+1}} y^{a_i} \nabla\phi_i\cdot \nabla \left(\frac{\zeta_i^2}{\phi_i} \right)d\bar x+ \int_{\partial\mathbb R_+^{n+1}} y^{a_i} \nabla\phi_i\cdot \nu \left(\frac{\zeta_i^2}{\phi_i} \right)dx =0.
\end{equation}
From this and the boundary term of equation (\ref{stabilityh}), we get 
 \begin{equation}\label{}
- 2\int_{\mathbb R_+^{n+1}} y^{a_i} \nabla\phi_i  \cdot \nabla\zeta_i \frac{\zeta_i}{\phi_i} d\bar x + \int_{\mathbb R_+^{n+1}} y^{a_i} |\nabla\phi_i|^2   \frac{\zeta_i^2}{\phi_i^2} d\bar x 
 + \int_{\partial\mathbb R_+^{n+1}} d_{s_i} \sum_{j=1}^m \partial_{j}H_i( v)\phi_j \frac{\zeta_i^2}{\phi_i}dx =0.
  \end{equation}
Applying the Young's inequality, for each index $i$, we obtain
 \begin{eqnarray}\label{intparHphi}
\int_{\partial\mathbb R_+^{n+1}} d_{s_i} \sum_{j=1}^m \partial_{j}H_i( v) \frac{\phi_j}{\phi_i} \zeta_i^2 d x&=&  \int_{\mathbb R_+^{n+1}} y^{a_i}  \left(   -  |\nabla\phi_i|^2   \frac{\zeta_i^2}{\phi_i^2} +2   \nabla\phi_i  \cdot \nabla\zeta_i \frac{\zeta_i}{\phi_i}  \right) d\bar x
\\&\le& \nonumber  \int_{\mathbb R_+^{n+1}} y^{a_i} |\nabla \zeta_i|^2 d\bar x   .
  \end{eqnarray}
 We now provide a lower bound for the  integrand in the left-hand side of (\ref{intparHphi}) when taking the sum on both sides of (\ref{intparHphi}) for $i=1,\cdots,m$, 
  \begin{eqnarray}
  \sum_{i,j} d_{s_i}  \partial_{j}H_i( v)  \frac{\phi_j}{\phi_i}\zeta_i^2 &=& 
  \sum_{i} d_{s_i} \partial_{i}H_i( v)  \zeta_i^2 + \sum_{i\neq j} d_{s_i}  \partial_{ j}H_i( v)   \frac{\phi_j}{\phi_i}\zeta_i^2
   \\&=& \nonumber \sum_{i} d_{s_i} \partial_{i}H_i( v)  \zeta_i^2 + \sum_{i < j} d_{s_i} \partial_{ j}H_i( v)   \frac{\phi_j}{\phi_i}\zeta_i^2 + \sum_{i> j} d_{s_i} \partial_{j}H_i( v)  \frac{\phi_j}{\phi_i}\zeta_i^2
    \\&= & \nonumber \sum_{i} d_{s_i} \partial_{i}H_i( v)  \zeta_i^2 + \sum_{i < j} d_{s_i} \partial_{j}H_i( v) \frac{\phi_j}{\phi_i}\zeta_i^2 + \sum_{i< j} d_{s_j} \partial_{i}H_j( v)  \frac{\phi_i}{\phi_j}\zeta_j^2 
    \\ &=& \nonumber \sum_{i} d_{s_i} \partial_{i}H_i( v)  \zeta_i^2 + \sum_{i < j}   (\phi_i \phi_j)^{-1} \left( d_{s_i} \partial_{j}H_i( v) \phi_j^2\zeta_i^2 +d_{s_j} \partial_{i}H_j( v) \phi_i^2\zeta_j^2\right) 
    \\&\ge &  \nonumber \sum_{i}  d_{s_i} \partial_{i}H_i ( v)  \zeta_i^2 + 2\sum_{i < j}   \sqrt{d_{s_i}d_{s_j}  \partial_{j}H_i( v) \partial_{i}H_j(v) } \zeta_i\zeta_j 
    \\&=& \nonumber \sum_{i,j}    \sqrt{d_{s_i}d_{s_j}  \partial_{j}H_i( v) \partial_{i}H_j(v) }   \zeta_i\zeta_j. 
  \end{eqnarray} 
This completes the proof. 

\hfill $ \Box$

%******   *******
 
 We now apply the stability inequality (\ref{stabilityinh}) to establish  a geometric Poincar\'{e} inequality.  Note that this inequality for the case of local scalar equations was first proved in \cite{sz} and it was used in \cite{fsv} and references therein to prove De Giorgi type results.  For the fractional Laplacian case this inequality was established in \cite{sv}.  We also refer interested readers to \cite{fg} and \cite{fs,dp} for the case of local and nonlocal systems, respectively.

\begin{thm}\label{lempoin}
 Assume that  $m,n\ge 1$ and $ v=(v_i)_i$ is a stable solution of (\ref{emainh}).  Then, for any $\eta=(\eta_k)_{k=1}^m \in C_c^1(\mathbb R_+^{n+1})$, the following inequality holds,
\begin{eqnarray}\label{poincare}
&&\sum_{i=1}^m  \frac{1}{d_{s_i}}   \int_{    \{|\nabla_{ x}  v_i|\neq 0 \} \cap \mathbb R_+^{n+1} }      y^{1-2s_i} \left(   |\nabla v_i|^2 \mathcal{A}_i^2 + | \nabla_{T_i} |\nabla_{ x}  v_i| |^2  \right)\eta_i^2 d\bar x
\\&&\nonumber +\sum_{i\neq j} \int_{\partial\mathbb R_+^{n+1}}   \left(  \sqrt{\partial_{j} H_i( v) \partial_{ i} H_j( v) }  |\nabla_{ x}  v_i|  |\nabla_{ x}  v_j| \eta_i \eta_j   - \partial_{j}H_i(v)  \nabla_{ x}  v_i \cdot   \nabla_{ x}  v_j \eta_i^2 \right) dx
\\&\le& \nonumber \sum_{i=1}^m \frac{1}{d_{s_i}}  \int_{\mathbb R_+^{n+1}} y^{1-2s_i}  |\nabla_{ x}  v_i|^2   |\nabla \eta_i|^2 d\bar x,
  \end{eqnarray} 
 where $\nabla_{T_i}$ stands for the tangential gradient along a given level set of $v_i$ and 
$\mathcal{A}_i^2$ for the sum of the squares of the principal curvatures of such a level set.
\end{thm}
\noindent\textbf{Proof:} Let  $v=(v_i)^m_{i=1}$ be a stable solution of (\ref{emainh}). From Lemma \ref{stabilitylem},  the stability inequality  (\ref{stabilityinh}) holds. Test the stability inequality with $\zeta_i:=|\nabla_{ x} v_i|\eta_i$ where each $\eta_i\in C^1_c(\mathbb R_+^{n+1})$, to get 
 \begin{eqnarray}\label{stablepoin}
I&:=& \sum_{i,j=1}^{m} \int_{\partial\mathbb{R}^{n+1}_+} \sqrt{d_{s_i}d_{s_j}  \partial_{j}H_i( v) \partial_{v_i}H_j(v) }   |\nabla_{ x} v_i|  |\nabla_{ x} v_j|\eta_i \eta_j  d x 
\\&\le& \nonumber \sum_{i=1}^{m} \int_{\mathbb {R}^{n+1}_+}  y^{a_i} |\nabla \left(   |\nabla_{ x} v_i|\eta_i   \right) |^2  d\bar x =:J.
  \end{eqnarray} 
Simplifying the left-hand side of the above inequality,  we get  
 \begin{eqnarray}\label{Ipoin}
I&=& \sum_{i=1}^{m} \int_{\partial\mathbb {R}^{n+1}_+} d_{s_i}   \partial_{v_i} H_i( v)  |\nabla_{ x} v_i|^2 \eta_i^2   d x  \\&& \nonumber +
\sum_{i\neq j} \int_{\partial\mathbb {R}^{n+1}_+} \sqrt{d_{s_i}d_{s_j}  \partial_{j}H_i( v) \partial_{v_i}H_j(v) }   |\nabla_{ x} v_i|  |\nabla_{ x} v_j|   \eta_i \eta_j  d x . 
  \end{eqnarray} 
Similarly,  for $J$ we have
\begin{equation}\label{Jpoin}
J=  \sum_{i=1}^{m} \int_{\mathbb {R}^{n+1}_+}  y^{a_i} \left[   |\nabla_{ x} v_i|^2 |\nabla \eta_i|^2 +    |\nabla|\nabla_{ x} v_i||^2  \eta_i^2  +\frac{1}{2}   \nabla |\nabla_{ x} v_i|^2  \cdot \nabla \eta_i^2  \right] d\bar x. 
  \end{equation} 
We now differentiate the $i^{th}$ equation of (\ref{emainh}) with respect to $x_k$ and multiply with $\partial_{x_k} v_i \eta_i^2$ for  $i=1,\cdots,m$ and $k=1,\cdots,n$, to get 
\begin{equation}
 \div(y^a_i   \nabla \partial_{x_k} v_i  ) \partial_{x_k} v_i \eta_i^2=0. 
  \end{equation} 
Integrating this by parts, we obtain
 \begin{eqnarray}\label{partsxk}
&& \int_{\mathbb {R}^{n+1}_+}  y^{a_i}  |\nabla \partial_{x_k} v_i|^2 \eta_i^2 d\bar x + \frac{1}{2} \int_{\mathbb{R}^{n+1}_+}  y^{a_i} \nabla |\partial_{x_k} v_i|^2 \cdot \nabla \eta_i^2 d\bar x
\\&=& \nonumber \int_{\partial\mathbb{R}^{n+1}_+} \lim_{y\to 0} y^{a_i}  (-\partial _y \partial_{x_k} v_i) \partial_{x_k} v_i \eta_i^2 dx. 
  \end{eqnarray} 
Differentiating the boundary term of (\ref{emainh}) with respect to $x_k$ yields 
\begin{equation}
 \lim_{y\to 0} y^{a_i}  (-\partial _y \partial_{x_k} v_i) \partial_{x_k} v_i \eta_i^2 = d_{s_i} \sum_{j=1}^m \partial_{j} H_i( v) \partial_{x_k} v_j \partial_{x_k} v_i   \eta_i^2.
   \end{equation} 
From this and (\ref{partsxk}),  we have
 \begin{eqnarray}\label{partsxkk}
&& \sum_{i=1}^{m} \int_{\mathbb {R}^{n+1}_+}  y^{a_i}  |\nabla \partial_{x_k} v_i|^2 \eta_i^2 d\bar x+  \sum_{i=1}^{m} \frac{1}{2} \int_{\mathbb{R}^{n+1}_+}  y^{a_i} \nabla |\partial_{x_k} v_i|^2 \cdot \nabla \eta_i^2  d\bar x
\\&=& \nonumber \sum_{i=1}^{m} d_{s_i} \int_{\partial\mathbb{R}^{n+1}_+}  \partial_{ v_i} H_i( v) |\partial_{x_k} v_i|^2   \eta_i^2  dx
 +  \sum_{i\neq j}^m  d_{s_i}  \int_{\partial\mathbb {R}^{n+1}_+}  \partial_{j} H_i( v) \partial_{x_k} v_j \partial_{x_k} v_i   \eta_i^2  dx. 
  \end{eqnarray} 
Taking sum on the index $k=1,\cdots, n$,  we get 
 \begin{eqnarray}\label{}
\ \ \ \ \ \  \sum_{i=1}^{m} d_{s_i} \int_{\partial\mathbb {R}^{n+1}_+}  \partial_{ v_i} H_i( v) |\nabla_{ x}  v_i|^2   \eta_i^2 dx &=&  \sum_{i=1}^{m} \int_{\mathbb {R}^{n+1}_+}  y^{a_i}   \sum_{k=1}^{n}   |\nabla \partial_{x_k} v_i|^2 \eta_i^2 d\bar x
  \\ && \nonumber + \frac{1}{2}  \sum_{i=1}^{m} \int_{\mathbb {R}^{n+1}_+}  y^{a_i} \nabla |\nabla_{ x} v_i|^2 \cdot \nabla \eta_i^2  d\bar x
  \\ &&  \nonumber -    \sum_{j \neq i}^m  d_{s_i}  \int_{\partial\mathbb {R}^{n+1}_+}  \partial_{j} H_i ( v)   \nabla_{ x} v_i \cdot \nabla_{ x} v_j  \eta_i^2  dx.
  \end{eqnarray} 
We now substitute this equality in $I\le J$ when  $I$ is given in (\ref{Ipoin}) and $J$ is given in (\ref{Jpoin}).  On can see that the term 
$ \frac{1}{2}  \sum_{i=1}^{m} \int_{\mathbb{R}^{n+1}_+}  y^{a_i} \nabla |\nabla_{ x} v_i|^2 \cdot \nabla \eta_i^2d\bar x  $    crosses  out and we end up with 
  \begin{eqnarray}\label{8left}
&&\sum_{i=1}^{m} \int_{\mathbb{R}^{n+1}_+ \cap \{ | \nabla v_i | \neq 0 \} }  y^{a_i}  \left( \sum_{k=1}^{n}   |\nabla \partial_{x_k} v_i|^2   -     |\nabla|\nabla_{ x} v_i||^2      \right) \eta_i^2 d\bar x
\\&& +  \sum_{i\neq j} \int_{\partial\mathbb {R}^{n+1}_+}     \left[  \sqrt{d_{s_i} d_{s_j} \partial_{ v_j} H_i( v) \partial_{ v_i} H_j( v) }   |\nabla_{ x} v_i|  |\nabla_{ x} v_j|   \eta_i \eta_j -  d_{s_i}  \partial_{j} H_i( v)   \nabla_{ x} v_i \cdot \nabla_{ x} v_j \eta_i^2  \right] dx
\\ &\le&  \sum_{i=1}^{m} \int_{\mathbb{R}^{n+1}_+ }  y^{a_i} |\nabla_{ x} v_i|^2 |\nabla \eta_i|^2d\bar x . 
 \end{eqnarray} 
According to formula (2.1) given in  \cite{sz}, the following geometric identity between the tangential gradients and curvatures holds. For any $w \in C^2(\Omega)$
 \begin{eqnarray}\label{8identity}
  \sum_{k=1}^{n} |\nabla \partial_k w|^2-|\nabla|\nabla w||^2=
\left\{
                      \begin{array}{ll}
                       |\nabla w|^2 (\sum_{l=1}^{n-1} \mathcal{\kappa}_l^2) +|\nabla_T|\nabla w||^2 & \hbox{for $x\in\{|\nabla w|>0\cap \Omega \}$,} \\
                       0 & \hbox{for $x\in\{|\nabla w|=0\cap \Omega \}$,}
                                                                       \end{array}
                    \right.   \end{eqnarray} 
 where $ \mathcal{\kappa}_l$ are the principal curvatures of the level set of $w$ at $ x$ and $\nabla_T$ denotes the orthogonal projection of the gradient along this level set. Applying  the above identity (\ref{8identity}) to  (\ref{8left}) completes the proof.
     
\hfill $ \Box$

%\begin{remark} Note that for symmetric systems  the boundary term in this inequality becomes 
%$$ \sum_{i\neq j} \int_{\partial\mathbb R_+^{n+1}}  \partial_{j} H_i( v)  \left(  |\nabla_{ x}  v_i|  |\nabla_{ x}  v_j| \eta_i \eta_j   -  \nabla_{ x}  v_i \cdot   \nabla_{ x}  v_j \eta_i^2 \right) dx$$
%Then by choosing an appropriate test function $\eta_i$ this term would have a fixed sign. 
%\end{remark}

%\subsection{Gradient and Energy Estimates} 

For the rest of this section,  we provide  energy and gradient estimates for solutions of (\ref{emainh}).  Then,  in next sections we apply these estimates  to establish De Giorgi type results and Liouville theorems.  Consider the  energy functional 
\begin{equation}\label{energy}
E_R(v):= \sum_{i=1}^{m}\frac{1}{2d_{s_i}} \int_{C_R}  y^{1-2s_i} |\nabla v_i|^2 d\bar x -  \int_{B_R} \tilde H( v)  dx ,
\end{equation}
when  $\partial_i \tilde H(v)=H_i(v)$ for $i=1,\cdots m$.    We finish this section by proving an energy estimate for $H$-monotone solutions of (\ref{emainh}).
\begin{thm}\label{energylayer}
Let  $ v=(v_i)_{i=1}^m$ be a bounded $H$-monotone solution of (\ref{emainh}) such that   
\begin{equation}
\lim_{x_n\to\infty} v_i( x',x_n,y)= L_i \ \ \ \text{for}\ \ \ x'\in\mathbb R^{n-1} \ \ \text{and}\ \ y\in\mathbb R^+,
\end{equation}
 when $L=(L_i)_{i=1}^m$ and each $L_i$ is a constant in $ \mathbb R$ and  $H(L)=0$.  Then,  the following energy estimates hold. 
\begin{enumerate}
\item[(i)] If $1/2<s_*<1$, then $E_R( v) \le  C R^{n-1}$,
\item[(ii)] If  $s_*=1/2$, then $E_R( v) \le  C R^{n-1}\log R$,
\item[(iii)] If $0<s_*<1/2$, then $E_R( v) \le  C R^{n-2s_*}$,
\end{enumerate} 
where the positive constant $C$ is independent from $R>1$. 
\end{thm}

 \noindent  \textbf{Proof:} Set the shift function $v_i^t( x,y):=v_i( x',x_n+t,y)$ for $( x',x_n,y)\in\mathbb R^{n+1}_+$ and $t\in\mathbb R$. It is straightforward to see that $v^t=(v^t_i)_{i=1}^m$ is a solution of (\ref{emainh}) and it satisfies pointwise estimates provided in Lemma \ref{asymp}. In addition,  for every parameter $t$ and  all indices $i$, one can see that $|v_i^t|\in {L^\infty}(\mathbb R^{n+1}_+)$. Therefore,  $v^t=(v^t_i)_{i=1}^m$  is a sequence of bounded functions and% for any $(x,y)\in \mathbb R^{n+1}_+$, 
  \begin{equation}
   \lim_{t\to\infty} |v_i^t( x,y)-L_i|+|\nabla v_i^t( x,y)|=0 \ \ \text{for all} \ \ (x,y)\in \mathbb R^{n+1}_+. 
   \end{equation}
This implies that 
  \begin{equation}
\lim_{t\to\infty} E_{R}( v^t)=0.
   \end{equation}
  We now differentiate the energy functional  in terms of parameter $t$ to get
 \begin{eqnarray}
\partial_t E_R( v^t)&=&  \sum_{i=1}^{m} \frac{1}{d_{s_i}}\int_{0}^{R}  \int_{B_R}  y^{1-2s_i} \nabla v^t_i\cdot\nabla(\partial_t v_i^t) d\bar x  -  \int_{B_R} \sum_{i=1}^{m}  H_{i}( v^t) \partial_t v_i^t dx
\\ &=& \nonumber -\sum_{i=1}^{m}\frac{1}{d_{s_i}}  \int_{0}^{R}  \int_{B_R} \div\left( y^{1-2s_i} \nabla v^t_i\right) \partial_t v_i^t d\bar x 
\\&& \nonumber + \sum_{i=1}^{m} \frac{1}{d_{s_i}} \int_{0}^{R}  \int_{\partial B_R} y^{1-2s_i} \left(\nabla v^t_i\cdot\nu\right) \partial_t v_i^t d\mathcal H^{n-1} dy
\\&&\nonumber  -  \sum_{i=1}^{m}  \int_{B_R} H_{i}( v^t) \partial_t v_i^t dx  . 
  \end{eqnarray}
 From the fact that $\div\left( y^{1-2s_i} \nabla v^t_i\right)=0$,  we can simplify the above as 
 \begin{equation}\label{}
 \partial_t E_R( v^t) =
 \sum_{i=1}^{m}\frac{1}{d_{s_i}}  \int_{0}^{R}  \int_{\partial B_R} y^{1-2s_i} \nabla v^t_i\cdot\nu \partial_t v_i^t d\mathcal H^{n-1} dy
 +  \sum_{i=1}^{m}\frac{1}{d_{s_i}}  \int_{B_R\times\{y=R\}}    y^{1-2s_i} \partial_y v^t_i \partial_t v_i^t   dx  .
   \end{equation}
Consider disjoint sets of indices $I $ and $J$ such that  $I\cup J=\{1,\cdots,m\}$ and   $ \partial_t v_\mu ^t >0> \partial_t v_\lambda ^t$ for $\mu\in I$ and $\lambda\in J$. Applying this,  we can expand  $\partial_t E_R( v^t)$ as
 \begin{eqnarray}\label{partERvt}
 \partial_t E_R( v^t)&=& \int_{\partial B_R}\int_{0}^{R}  \left( \sum_{\mu\in I} \frac{1}{d_{s_\mu}} y^{1-2s_\mu} \partial_\nu v^t_\mu  \partial_t v_\mu^t + \sum_{\lambda\in J} \frac{1}{d_{s_\lambda}} y^{1-2s_\lambda} \partial_\nu v^t_\lambda  \partial_t v_\lambda^t  \right) dy d\mathcal H^{n-1}
 \\&& \nonumber + \int_{B_R\times\{y=R\}} \left( \sum_{\mu\in I} \frac{1}{d_{s_\mu}} y^{1-2s_\mu} \partial_y v^t_\mu  \partial_t v_\mu^t + \sum_{\lambda\in J}  \frac{1}{d_{s_\lambda}} y^{1-2s_\lambda} \partial_y v^t_\lambda  \partial_t v_\lambda^t  \right) dx . 
   \end{eqnarray}
Lemma \ref{asymp} implies that there exists a constant $M$ such that $ |\partial_\nu v^t_i | \le \frac{M}{1+y}$ for $(x,y)\in\mathbb{R}^{n+1}_+$ and $ |\partial_y v_i^t| \le \frac{M}{y}$  for  $x\in\mathbb{R}^{n}$ and $ y>1$. From these estimates and (\ref{partERvt}),   we get  
\begin{eqnarray}\label{tEt}
 \partial_t E_R( v^t)&\ge& M \int_{\partial B_R}\int_{0}^{R} \sum_{\mu\in I}  \frac{1}{d_{s_\mu}} y^{1-2s_\mu} \left(- \frac{1}{1+y} \right)  \partial_t v_\mu^t  dy d\mathcal H^{n-1}
 \\&&\nonumber +  M \int_{\partial B_R}\int_{0}^{R}  \sum_{\lambda\in J} \frac{1}{d_{s_\lambda}}  y^{1-2s_\lambda} \left( \frac{1}{1+y} \right)  \partial_t v_\lambda^t   dy  d\mathcal H^{n-1}
 \\&&\nonumber + M \int_{B_R\times\{y=R\}}  \sum_{\mu\in I}  \frac{1}{d_{s_\mu}}  y^{1-2s_\mu} \left(- \frac{1}{y} \right)  \partial_t v_\mu^t  dx
 \\&&\nonumber + M \int_{B_R\times\{y=R\}} \sum_{\lambda\in J} \frac{1}{d_{s_\lambda}}  y^{1-2s_\lambda} \left(\frac{1}{y} \right)  \partial_t v_\lambda^t   dx . 
   \end{eqnarray}
Note that  $ E_R( v)= E_R( v^T)-\int_0^T \partial_t E_R( v^t) dt$ for every $T>0$.  Combining  this and (\ref{tEt}), we conclude  
   %Integration of the energy with respect to time implies that for every $T>0$,  we have $ E_R( v)= E_R( v^T)-\int_0^T \partial_t E_R( v^t) dt$.   Therefore,   
\begin{eqnarray}\label{ERvine}
E_R( v) &\le& E_R( v^T)  +M  \int_{\partial B_R}\int_{0}^{R}   \sum_{\mu\in I} \frac{1}{d_{s_\mu}} \left(\frac{y^{1-2s_\mu} }{1+y} \right)  \int_0^T \partial_t v_\mu^t dt dy d\mathcal H^{n-1}
\\&&\nonumber - M  \int_{\partial B_R}\int_{0}^{R}  \sum_{\lambda\in J} \frac{1}{d_{s_\lambda}} \left( \frac{y^{1-2s_\lambda} }{1+y} \right) \int_0^T \partial_t v_\lambda^t dt  dy d\mathcal H^{n-1}
 \\&&\nonumber + M \int_{B_R\times\{y=R\}}  \sum_{\mu\in I}\frac{1}{d_{s_\mu}} y^{1-2s_\mu-1}   \int_0^T \partial_t v_\mu^t dt  dx
 \\&&\nonumber  - M   \int_{B_R\times\{y=R\}}\sum_{\lambda\in J}  \frac{1}{d_{s_\lambda}} y^{1-2s_\lambda-1}    \int_0^T \partial_t v_\lambda^t  dt  dx.
   \end{eqnarray}
Simplifying (\ref{ERvine}) yields  
\begin{eqnarray}
E_R( v) &\le& E_R( v^T) + M \int_{\partial B_R} \int_{0}^{R} \sum_{\mu\in I}  \frac{y^{1-2s_\mu} }{1+y}      (v_\mu^T - v_\mu)  dy d\mathcal H^{n-1}
\\ && \nonumber +  M \int_{\partial B_R} \int_{0}^{R} \sum_{\lambda\in J}  \frac{y^{1-2s_\lambda} }{1+y}   (v_\lambda-v_\lambda^T) dy d\mathcal H^{n-1}
\\&& \nonumber + M \int_{B_R\times\{y=R\}}  \sum_{\mu\in I} y^{-2s_\mu}   (v_\mu^T -v_\mu)  dx 
 +  M \int_{B_R\times\{y=R\}}  \sum_{\lambda\in J} y^{-2s_\lambda}   (v_\lambda-v_\lambda^T) dx. 
   \end{eqnarray}
Note that for  $\mu\in I$ and $\lambda\in J$,  we  have  $v_\mu^T \ge v_\mu$ and $v_\lambda \ge v_\lambda^T$.  Therefore, 
 \begin{equation}
E_R(v) \le  E_R(v^T) + M \int_{\partial B_R}  \int_0^R  \sum_{i=1}^{m} \frac{y^{1-2s_i} }{1+y} dy d\mathcal H^{n-1} + M  \int_{B_R\times\{y=R\}}  \sum_{i=1}^{m} y^{-2s_i} dx. 
   \end{equation}
Let  $T\to\infty$,  then $E_R(v^T)$ approaches zero.  Doing integration by parts we obtain 
   \begin{equation}
E_R( v) \le  M  \sum_{i=1}^{m} \left[ R^{n-2s_i} \chi_{\{0<s_i<1/2\}} +  R^{n-1} \chi_{\{1/2<s_i<1\}} 
 +   R^{n-1}\log R \chi_{\{s_i=1/2\}}   +   R^{n-2s_i} \right] . 
   \end{equation}   
This completes the proof.

               \hfill $ \Box$

Note that ideas and techniques applied in the above proof are strongly motivated by the ones provided by Ambrosia and Cabr\'{e} in \cite{ac} and by Cabr\'{e} and Cinti in \cite{ccinti,ccinti1} for local and nonlocal scalar equations, respectively. We would like to refer interested readers to  \cite{fg} by Ghoussoub and the author and  to \cite{fs} by Sire and the author for local and nonlocal system of equations, respectively.  In the next two lemmata, we prove gradient estimates for bounded solutions of (\ref{emainh}).

%For the next lemma the only assumption on solutions is boundedness and the nonlinearity $H$ can be any function with sufficient regularity.   We shall apply this in next sections to prove Liouville theorems.     

\begin{lemma}\label{deltavx} 
Suppose that $v=(v_i)_{i=1}^m$ is a bounded solution of (\ref{emainh}). Then, the following gradient estimate holds for any $i=1,\cdots, m$  
\begin{equation}
\int_{C_R} y^{a_i} |\nabla v_i|^2d\bar x  \le  C R^n,
\end{equation}
when   $C$ is a positive constant independent from $R$. 
\end{lemma}
 \noindent\textbf{Proof:} Let $v$ be a bounded solution of   (\ref{emainh}). Multiplying (\ref{emainh})   by $v_i\psi_R$ where $\psi_R$ is a test function and integrating over $\mathbb R^{n+1}_+$, we get  
\begin{equation}
  \int_{\mathbb R^{n+1}_+} y^{a_i} |\nabla v_i|^2 \psi_R^2 d\bar x + 2   \int_{\mathbb R^{n+1}_+} y^{a_i} \nabla \psi_R\cdot\nabla v_i \psi_R v_i d\bar x =   d_{s_i}\int_{\partial\mathbb R^{n+1}_+}  H_i(v) v_i \psi_R^2 d x.
  \end{equation}
Using the Cauchy-Schwarz inequality,  we conclude
\begin{equation}
 2 \int_{\mathbb R^{n+1}_+} y^{a_i} \psi_R \nabla \psi_R\cdot\nabla v_i  v_i d\bar x \le \frac{1}{2} \int_{\mathbb R^{n+1}_{+}} y^{a_i} |\nabla v_i|^2 \psi_R^2 d\bar x + 2 \int_{\mathbb R^{n+1}_+} y^{a_i} |\nabla \psi_R|^2 v_i^2 d\bar x .\end{equation}
From this,  we obtain
\begin{equation}
  \int_{\mathbb R^{n+1}_+} y^{a_i} |\nabla v_i|^2 \psi_R^2 d\bar x \le 4  \int_{\mathbb R^{n+1}_+} y^{a_i} |\nabla \psi_R|^2 v_i^2d\bar x + 2  d_{s_i}  \int_{\partial  \mathbb R^{n+1}_+}  H_i( v) v_i \psi_R^2 d x.\end{equation}
From the boundedness of $v$ and the fact that $H_i\in C^{1,\gamma}$, we obtain the following bound by choosing an  appropriate test function, 
\begin{equation}
 \int_{C_R} y^{a_i} |\nabla v_i|^2 \psi_R^2  d\bar x \le C R^{n-2s_i} + CR^n \le C R^n, \end{equation}
where the constant $C$ is independent from $R$ but may depend on $v=(v_i)_{i=1}^m$. This completes the proof.

  \hfill $ \Box$
  
  We now assume extra conditions on the sign of the nonlinearity $H$. This enables us to  prove stronger gradient estimates on solutions of (\ref{emainh}).  
    
  \begin{lemma}\label{Hiv}
Suppose that $ v=(v_i)_{i=1}^m$ is a bounded solution of (\ref{emainh}). 
\begin{enumerate}
\item[(i)] If  $H_i(v) \ge0$ for any $1\le i\le m$, then  
%\begin{equation}\label{Hpositive}
$\int_{C_R} y^{1-2s_i} |\nabla v_i|^2 d\bar x \le C R^{n-2s_i}$. 
%\end{equation}  
\item[(ii)] If  $\sum_{i=1}^m v_i H_i(v)\le 0$, then %the following gradient estimate  holds 
%\begin{equation}\label{vHpositive}
$\sum_{i=1}^m \int_{C_R} y^{1-2s_i} |\nabla v_i|^2 d\bar x \le  C R^{n-2s_*}.$ 
%\end{equation}
\end{enumerate} 
Here, the positive constant $C$  is independent from $R$. 
\end{lemma}
 \noindent\textbf{Proof:} We start the proof with Part (i). Suppose that for a particular index $1\le i\le m$ we have $H_i(v) \ge0$. We now multiply  (\ref{emainh}) with $(v_i-||v_i||_{\infty})  \psi_R$ when  $\psi_R$ is a positive test function. % and integrating,  we get  
%\begin{equation}\label{EnHpos} \int_{\mathbb R^{n+1}_+}  (v_i-||v_i||_{\infty})  \psi_R \div(y^{a_i} \nabla v_i)d\bar x =0. \end{equation}
Applying integration by parts and using the fact that $H_i$ has a fixed sign,  we obtain 
\begin{eqnarray}\label{EnHpos1}
\int_{\mathbb R^{n+1}_+} y^{a_i}  \nabla[ (v_i-||v_i||_{\infty})  \psi_R]\cdot\nabla v_i d\bar x  &=&\int_{\partial\mathbb R^{n+1}_+}  \lim_{y\to 0} y^{a_i} (-\partial_y v_i)  (v_i-||v_i||_{\infty})   \psi_R d x
\\&=&\nonumber d_{s_i}\int_{\partial\mathbb R^{n+1}_+}   (v_i-||v_i||_{\infty})  H_i(v) \psi_R d x
\le 0.
\end{eqnarray}
This implies that 
\begin{eqnarray}
\int_{\mathbb R^{n+1}_+} y^{a_i} |\nabla v_i|^2 \psi_R d\bar x &\le& - \int_{\mathbb R^{n+1}_+} y^{a_i} (v_i-||v_i||_{\infty})   \nabla \psi_R\cdot \nabla v_i d\bar x
\\&=&\nonumber  \frac{1}{2}  \int_{\mathbb R^{n+1}_+} y^{a_i} \nabla (v_i-||v_i||_{\infty})^2\cdot \nabla \psi_R  d\bar x
\\&=& \nonumber - \frac{1}{2}  \int_{\mathbb R^{n+1}_+} (v_i-||v_i||_{\infty})^2 \div( y^{a_i}\nabla \psi_R) d\bar x
\\&& \nonumber + \frac{1}{2} \int_{\partial\mathbb R^{n+1}_+}  (-\lim_{y\to 0}  y^{a_i} \partial_y \psi_R) (v_i-||v_i||_{\infty})^2d x
\\&=:&\nonumber  I+J . 
\end{eqnarray}
From boundedness of  $v_i$ and applying an appropriate test function $\psi_R$, we obtain  
\begin{eqnarray}\label{Ivpsi}
I&=&- \frac{1}{2}  \int_{\mathbb R^{n+1}_+} (v_i-||v_i||_{\infty})^2 \left( y^{a_i} \Delta_{x} \psi_R +a_i y^{a_i-1} \partial_y \psi_R   \right) d\bar x
\\&\le& \nonumber C  \int_{C_{2R}} \left(\frac{y^{a_i}}{R^2} +\frac{y^{a_i-1}}{R} \right)  d \bar x \le C R^{n} \int_{R}^{2R}  \left(\frac{y^{a_i}}{R^2} +\frac{y^{a_i-1}}{R} \right)  dy \le \nonumber C R^{n-2s_i} , 
\end{eqnarray}
when  the positive constant $C$ is independent from $R$. In  above,  we have used the fact that the test function $\psi_R$ can be chosen such that $|\partial_y \psi_R |\le CR^{-1}$ and $|\Delta_{ x} \psi_R| \le C R^{-2}$ in $C_R$. Using these properties,  we conclude  $ |y^{a_i}  \partial_y \psi_R |\le C R^{-2s_i}$ in $C_R$.  Similarly, applying the test function $\psi_R$ with the above properties, we obtain the following upper bound for $J$ given in  (\ref{Ivpsi}) 
\begin{equation}
 J \le C R^{n-2s_i}.
 \end{equation}
This completes the proof of Part (i). We now provide a proof for Part (ii). %Now suppose that   $\sum_{i=1}^m v_i H_i(v)\ge 0$. 
Multiplying the $i^{th}$ equation of (\ref{emainh}) with $v_i \psi_R$, when  $\psi_R$ is the same test function applied in Part (i), and integrating we get  
\begin{eqnarray}\label{EnHpos2}
 \ \ \ \ \ \sum_{i=1}^m \frac{1}{d_{s_i}}\int_{\mathbb R^{n+1}_+} y^{a_i}  \nabla[ v_i  \psi_R]\cdot\nabla v_i d \bar x  &=& \sum_{i=1}^m \int_{\partial\mathbb R^{n+1}_+}  \lim_{y\to 0} y^{a_i} (-\partial_y v_i)  v_i   \psi_R dx
\\&=&  \nonumber \int_{\partial\mathbb R^{n+1}_+}   \sum_{i=1}^m v_i  H_i(v) \psi_R dx   .
\end{eqnarray}
Note that the latter is nonpositive since $\sum_{i=1}^m v_i H_i(v)\le 0$.  The rest of the proof is similar to Part (i) and we omit it here.  
 
               \hfill $ \Box$
               
 We end this section with a monotonicity formula for bounded solutions of \eqref{emainh}. Consider the following function for $R\ge 1$,   
 \begin{equation}\label{monIR}
I(R):= R^{2s^*-n} \left(\frac{1}{2} \int_{B_R^+}  \sum_{i=1}^{m} y^{a_i} |\nabla v_i|^2 d x dy - \int_{B_R\times\{y=0\}} \tilde H( v) d x \right). 
\end{equation}
when $\partial_i \tilde H(v)= H_i(v)$ for every $i$.  We show that $I(R)$ is a nondecreasing function of $R$,   when   $\tilde H \le 0$,  that is equivalent to  
 \begin{equation}\label{IRI1} 
I(R)\ge I(1) \ \ \text{for} \ \  R>1.
 \end{equation}
 Note that from definitions of $I(R)$ in (\ref{monIR}) and $E_R(v)$ in (\ref{energy}), we have
 \begin{equation} \label{Ervrn}
 E_R(v) \ge R^{n-2s^*} I(R). 
 \end{equation}
 From (\ref{Ervrn}) and (\ref{IRI1}), we conclude 
 \begin{equation} 
\sum_{i=1}^{m}\frac{1}{2d_{s_i}} \int_{C_R}  y^{1-2s_i} |\nabla v_i|^2 d\bar x -  \int_{B_R} \tilde H( v)  dx  \ge C R^{n-2s^*},
 \end{equation}
and the positive constant $C$ is independent from $R$.  This  clarifies the sharpness of gradient and energy estimates provided in this section, in particular when $0<s_*=s^*<1/2$.

 \begin{thm}\label{thmmonoI} Let $v=(v_i)_{i=1}^m$ be a bounded solution of \eqref{emainh}. Then, $I(R)$ 
%\begin{equation}
%I(R)= R^{2s^*-n} \left(\frac{1}{2} \int_{B_R^+}  \sum_{i=1}^{m} y^{a_i} |\nabla v_i|^2 d x dy - \int_{B_R\times\{y=0\}} \tilde H( v) d x \right) \end{equation}
 is a nondecreasing function of $R\ge 1$ when $\partial_i \tilde H(v)= H_i(v)$ for every $1\le i\le m$ and  $\tilde H \le 0$.    
\end{thm}                
 \noindent  \textbf{Proof:}     
Differentiating $I(R)$ with respect to $R$,  yields 
 \begin{eqnarray}\label{IprimeR}
I'(R) &=& - \frac{n-2s^*}{2} R^{2s^*-n-1}  \int_{B_R^+} \sum_{i=1}^{m}  y^{a_i} |\nabla v_i|^2 d x dy + (n-2s^*)  R^{2s^*-n-1} \int_{B_R \times \{y=0\}} \tilde H(  v) d x 
\\&&+ \nonumber
\frac{1}{2}    R^{2s^*-n}  \int_{\partial^+{B_R^+}} \sum_{i=1}^{m}  y^{a_i} |\nabla v_i|^2 d x dy - R^{2s^*-n} \int_{\partial {B_R}\times\{y=0\}} \tilde H( v) d x  .
    \end{eqnarray}
 Applying some standard arguments in regards to the Pohozaev identity implies that 
  \begin{eqnarray}\label{n2nablav}
 \frac{n}{2} \int_{B_R^+}   \sum_{i=1}^{m} y^{a_i} |\nabla v_i|^2  &=&   \sum_{i=1}^{m} s_i \int_{B_R^+}     y^{a_i} |\nabla v_i|^2  + n \int_{B_R} \tilde H( v) d x - R \int_{\partial B_R} \tilde H( v) d \mathcal{H}^{n-1}  
 \\&& \nonumber - R  \int_{\partial^+ B_R^+}  \sum_{i=1}^{m} y^{a_i} (\partial_{ \nu} v_i )^2 d \mathcal{H}^{n} +  
 \frac{R}{2} \int_{\partial^+{B_R^+}} \sum_{i=1}^{m} y^{a_i}   |\nabla v_i|^2 .
   \end{eqnarray}
Combining (\ref{IprimeR}) and (\ref{n2nablav}), we conclude 
\begin{eqnarray}
I'(R) R^{n+1-2 s^*} &=& R  \int_{\partial^+ B_R^+}  \sum_{i=1}^{m} y^{a_i} (\partial_\nu v_i )^2 d \mathcal{H}^{n} 
\\&&\nonumber +   \int_{B_R^+}  \sum_{i=1}^{m} (s^*-s_i)  y^{a_i} |\nabla v_i|^2 - 2 s^* \int_{B_R\times\{y=0\}} \tilde H( v) d x.
   \end{eqnarray}
Since $\tilde H\le 0$ and $s^*\ge s_i$ for every $1\le i\le m$,  we  have $I'(R)\ge 0$.

 \hfill $ \Box$          
               
   % ****Alikakos for Liouville  theorem....                   

 %%%%%%%%%%%%%%%%%%%%%%%%%
 %%%%%%%%%%%%%%%%%%%%%%%%%%

\section{Symmetry Results and Liouville Theorems; General Nonlinearity}\label{nonlin}

In this section, we provide De Giorgi type results for stable and $H$-monotone solutions of symmetric system (\ref{emainh}) with a general nonlinearity in lower dimensions.   Just like in the proof of the classical De Giorgi's conjecture in dimensions $n=2,3$  providing  a Liouville theorem for the quotient of partial derivatives  is a milestone.  To be mathematically more precise,  consider the case of scalar equation (\ref{sceqH}), it was observed  by Berestycki, Caffarelli and Nirenberg in \cite{bcn},  by Ghoussoub and Gui in \cite{gg} and  by Ambrosio and Cabr\'{e} in \cite{ac} that if  $\phi\in L_{loc}(\mathbb R^n)$ such that $\phi^2>0$ a.e. and $\sigma \in H^1_{loc}(\mathbb R^n)$ satisfy 
\begin{equation}\label{divphi} -\sigma \div(\phi^2 \sigma) \le 0 \ \ \ \ \text{in }\ \ \mathbb R^n, \end{equation} 
under the decay-growth assumption 
\begin{equation}\label{intsigmaR} 
\int_{B_{2R}\setminus B_R} \phi^2 \sigma^2 < C R^2\ \ \ \ \text{for} \ \ R>1,
\end{equation}
 then $\sigma$ must be a constant.  Eventually, $\sigma $ is set to be $\sigma=\frac{\nabla u \cdot \eta}{\partial_n u}$ for an arbitrary direction $\eta\in\mathbb R^n$ and $\phi=\partial_n u$.   Cabr\'{e} and Sire in \cite{cs1,cs2} proved a similar Liouville theorem for nonlocal scalar equations  that is (\ref{mainh}) for $m=1$.  Most recently, Ghousssoub and the author in \cite{fg} and later Sire and the author in \cite{fs}  provided counterparts of this Liouville theorem for the local and nonlocal gradient system (\ref{syeqH}), respectively. 
 % We now provide a nonlocal Liouville theorem,  given in \cite{fs}, that we apply to prove our main results. 
  %A nonlocal version of this, for the case of $m=1$, is given by Cabr\'{e} and Sire in \cite{cs1,cs2}. 
 %For the case of systems (\ref{emainh}) providing such a Liouville theorem is slightly tricky in the sense that one needs to figure a counterpart of (\ref{divphi}) to begin with. In other words, only considering  $-\sigma_i \div(\phi_i^2 \sigma_i) \le 0$ for $i=1,\cdots,m$ does not help since this is a decoupled system. We provide the needed Liouville theorem as Theorem \ref{liouville} and we apply it to establish De Giorgi type results and Liouville theorem for (\ref{emainh}).  Note that the concept of symmetric systems seems to be crucial when we apply this Liouville theorem to prove De Giorgi type results for system (\ref{mainh}).   
  Consider the following set of functions $\mathcal F$ that  is somewhat standard in this  context 
   $$\mathcal F:=\left\{F:\mathbb R^+\to\mathbb R^+, F \  \text{is nondecreasing and} \ \int_{2}^{\infty} \frac{dr}{rF(r)}=\infty\right\}.$$
For example, $F(r)=\ln r$ and $F(r)\equiv constant>0$ belong to this class and $F(r)=r$ does not belong to $\mathcal F$. As far as we know,  this set of functions was introduced by Karp in \cite{k1,k2} and used in \cite{mos,fs,mf,ccinti1,ccinti,cs1,cs2}. 
 
 \begin{thm}\label{liouville} Assume that   $\phi_i \in L^{\infty}_{loc} (\overline{\mathbb{R}_+^{n+1}}) $ is a positive function and $\sigma_i \in H^1_{loc}(\overline{\mathbb{R}_+^{n+1}},y^{a_i}) $ satisfies
\begin{equation}\label{liouassum}
 \limsup_{R\to\infty} \frac{1}{R^2F(R)} \int_{C_R}  \sum_{i=1}^{m}  y^{a_i}  \phi_i^2\sigma_i^2 d\bar x <\infty,
 \end{equation}
for $F\in \mathcal F$ and    $i=1,\cdots,m$.  If the sequence $(\sigma_i)_{i=1}^m$ is a solution of 
 \begin{eqnarray}\label{div}
 \left\{ \begin{array}{lcl}
\hfill  -\sigma_i \div(y^{a_i}\phi_i^2\nabla \sigma_i)&\le& 0   \ \ \text{in}\ \ \mathbb{R}_+^{n+1},\\   
\hfill -\lim_{y\to0}y^{a_i} \phi_i^2 \sigma_i \partial_{y} \sigma_i &\le&  \sum_{j=1}^{n} h_{ij} f(\sigma_j-\sigma_i)\sigma_i  \ \ \text{in}\ \ \partial\mathbb{R}_+^{n+1},
\end{array}\right.
  \end{eqnarray}
when $0\le h_{i,j}\in L_{loc}^1(\mathbb{R}^n)$, $h_{i,j}=h_{j,i}$ and $f\in L_{loc}^1(\mathbb{R})$ is an odd function such that $f(t)\ge 0$ for $t\in\mathbb R^+$.  Then,  each function $\sigma_i$ is constant for all  $i=1,\cdots,m$.
\end{thm}

In this article, we apply the above theorem frequently to prove our main results for solutions of (\ref{emainh}).    
 Suppose that $v=(v_i)_{i=1}^m$ is a $H$-monotone solution of (\ref{emainh}). 
%In this section,  we study (\ref{emainh}) with a general nonlinearity.   
Let  $\phi_i := \partial _n v_i$ and $\psi_i:=\nabla v_i\cdot\eta$ for any fixed $\eta=(\eta',0)\in \mathbb{R}^{n-1}\times\{0\}$. Then, sequences $(w_i)_{i=1}^m=(\phi_i)_{i=1}^m$ and $(w_i)_{i=1}^m=(\psi_i)_{i=1}^m$ satisfy the following linearized equation
\begin{eqnarray}\label{pointstability1}
 \left\{ \begin{array}{lcl}
\hfill \div(y^{a_i} \nabla w_i)&=& 0   \ \ \text{in}\ \ \mathbb{R}_+^{n+1},\\   
\hfill -\lim_{y\to0}y^{a_i} \partial_{y} w_i&=& d_{s_i} \sum_{j=1}^m \partial_{j}H_i( v)w_j   \ \ \text{in}\ \ \partial\mathbb{R}_+^{n+1} .
\end{array}\right.
  \end{eqnarray}
 Since $v$ is a $H$-monotone solution,  $\phi_i$ does not change sign for each $i$. Then,  the sequence of functions  $\sigma=(\sigma_i)_{i=1}^m$ for $\sigma_i:=\frac{\psi_i}{\phi_i}$ satisfies  the following 
  \begin{equation}\label{sigma1}
 \div(y^{a_i} \phi_i^2 \nabla \sigma_i) = \div(y^{a_i} [\phi_i \nabla \sigma_i-\sigma_i \nabla \phi_i])
 = \phi_i  \div(y^{a_i}\nabla \psi_i) - \psi_i  \div(y^{a_i}\nabla \phi_i) 
 =  0. 
  \end{equation}
On the other hand,  in $\partial\mathbb {R}_+^{n+1}$ we have 
\begin{equation}\label{}
   -\lim_{y\to0}y^{a_i} \partial_{y} \sigma_i=   -\lim_{y\to0}y^{a_i} \partial_{y} \psi_i  \phi_i^{-1} + \lim_{y\to0}y^{a_i} \partial_{y} \phi_i \frac{\psi_i}{\phi_i^2} 
   =   d_{s_i} \sum_{j=1}^m \partial_{j}H_i( v) \left[  \frac{\psi_j}{\phi_i} - \phi_j \frac{\psi_i}{\phi_i^2} \right], 
     \end{equation}
     which implies that   
\begin{equation}\label{sigma2}
     -\lim_{y\to0}y^{a_i} \phi_i^2 \sigma_i \partial_{y} \sigma_i = d_{s_i} \sum_{j=1}^m \partial_{j}H_i( v) [  \psi_j \phi_i  - \phi_j  \psi_i ]\sigma_i
     = 
      d_{s_i} \sum_{j=1}^m \partial_{j}H_i( v) \phi_i \phi_j   (\sigma_j-\sigma_i)\sigma_i   . \end{equation}
Note that if we set $h_{i,j}=\partial_{j}H_i( v) \phi_i \phi_j $ and $f$ to be  the identity function, then the above equations (\ref{sigma1}) and (\ref{sigma2}) satisfy (\ref{div}). 
%This gives us the chance to apply the following Liouville theorem given by Sire and the author  in \cite{fs} when $f$ is the identity function and $h_{i,j}=\partial_{j}H_i( v) \phi_i \phi_j $. 
Note that for symmetric systems we have $h_{i,j}=h_{j,i}$ and for $H$-monotone solutions we have $h_{i,j}>0$.   %Therefore, for symmetric systems the following  quantity has a fixed sing. 
The following computation for the right-hand side of (\ref{sigma2}) is our main observation to define symmetric systems and $H$-monotone solutions and to establish Theorem \ref{liouville};
\begin{eqnarray}
 \sum_{i,j=1}^m \phi_i \phi_j \partial_j H_i(v) \sigma_i (\sigma_i-\sigma_j)&=& \sum_{i< j}   \phi_i \phi_j \partial_j H_i(v)  \sigma_i (\sigma_i-\sigma_j) + \sum_{i> j}  \phi_i \phi_j \partial_j H_i(v)  \sigma_i  (\sigma_i-\sigma_j)  
 \\&=& \nonumber \sum_{i< j}  \phi_i \phi_j \partial_j H_i(u)   \sigma_i (\sigma_i-\sigma_j) + \sum_{i< j}  \phi_i \phi_j \partial_j H_i(v)  \sigma_j (\sigma_j-\sigma_i)  
 \\&=&\nonumber  \sum_{i< j}  \phi_i \phi_j \partial_j H_i(v)  (\sigma_i-\sigma_j)^2 \ge 0. 
\end{eqnarray}
We are now ready to state the following De Giorgi type result. 

\begin{thm}\label{thsymv}
Suppose that $ v=(v_i)_{i=1}^m$ is a bounded solution of  the orientable symmetric system \eqref{emainh}. Assume also that   either $n=2$, $0<s_i<1$ and  $ v$ is stable or $n = 3$,  $1/2 \le s_* <1$ and $ v$ is $H$-monotone. Then, there exist a constant $ \Gamma_i\in{S}^{n-1}$ and $v^*_i: \mathbb{R^+}\times  \mathbb{R}^+\to  \mathbb{R}$ such that
  \begin{equation}
   v_i( x,y)=v^*_i(  \Gamma_i\cdot x,y) \ \ \text{for} \ \ (x,y)\in\mathbb R^{n+1}_+,
    \end{equation}
 and $i=1,\cdots,m$. Moreover, for all $i,j=1,\cdots, m$ vectors $\nabla_x v_i(x,0)$ and $\nabla_x v_j (x,0)$ are parallel and the angle between two vectors is $\arccos\left(\frac{|\partial_j H_i(v)|}{\partial_j H_i(v)}\right)$. 
\end{thm}

\noindent\textbf{Proof:}  First,  let $n=2$ and $v$ be a stable solution of (\ref{emainh}).   We apply Lemma  \ref{deltavx} and Theorem \ref{liouville} when $h_{i,j}= \partial_{j} H_i(v) \partial_n v_i \partial_n v_j$, $f$ is the identity function and $F$ is constant.  % Note that this result holds for the weaker assumption of stability on solutions as well, due to the geometric Poincar\'{e} inequality given in Theorem \ref{lempoin} and  Lemma \ref{deltavx}.
 For each $1\le k\le m$,  set  $\eta_k (\bar x) :=\rho_R (\bar x) $  in  the geometric Poincar\'{e} inequality (\ref{poincare}) 
%sign(\partial_{n} v_k) 
 when  
\begin{equation}\label{testchi}
\rho_R (\bar x):=\left\{
                      \begin{array}{ll}
                       \log R, & \hbox{if $|\bar x|\le\sqrt{R}$,} \\
                     2 \frac{ \log R - \log |\bar x|}{\log R}, & \hbox{if $\sqrt{R}< |\bar x|< R$,} \\
                       0, & \hbox{if $|\bar x|\ge R$.}
                                                                       \end{array}
                    \right.
                  \end{equation}  
From (\ref{poincare}) and the fact that $|\nabla \rho_R|\le \frac{C}{|\bar x|}$ when $ \sqrt{R} <\bar x < R$,  %for a positive constant $C$ that is independent from $R$.  
we conclude 
\begin{eqnarray}\label{inelogrr}
&& |\log R|^2 \sum_{i=1}^m  \frac{1}{d_{s_i}}   \int_{  B^+_{\sqrt R} \cap \mathbb R_+^{n+1} }      y^{1-2s_i} \left(   |\nabla v_i|^2 \mathcal{A}_i^2 + | \nabla_{T_i} |\nabla_{ x}  v_i| |^2  \right) d\bar x
\\&& \nonumber +|\log R|^2  \sum_{i\neq j} \int_{ B^+_{\sqrt R} \cap   \partial\mathbb R_+^{n+1}}   \left(  |\partial_{j} H_i( v)|  |\nabla_{ x}  v_i|  |\nabla_{ x}  v_j|    - \partial_{j}H_i(v)  \nabla_{ x}  v_i \cdot   \nabla_{ x}  v_j  \right) dx
\\&\le& \nonumber C \int_{B^+_{R}\setminus B^+_{\sqrt R}}  \frac{y^{1-2s_i} |\nabla v_i|^2}{|\bar x|^2} d \bar x . 
\end{eqnarray}
Here,  we have used the notion of symmetric systems that is $\sqrt{\partial_{j} H_i( v) \partial_{ i} H_j( v)} =|\partial_{ i} H_j( v)|$. On the other hand, Lemma \ref{deltavx} implies that  for each index $1\le i\le m$, 
\begin{equation}\label{mm}
\int_{B_R^+} y^{1-2s_i} |\nabla v_i|^2 d\bar x \le   C R^2. 
\end{equation}
Straightforward calculations show that for each $i$,  we have  
\begin{eqnarray}
&&\int_{B^+_{R}\setminus B^+_{\sqrt R}}  \frac{1}{ |\bar x|^2}y^{1-2s_i} |\nabla v_i|^2 d \bar x
\\&=& \nonumber 2\int_{B^+_{R}\setminus B^+_{\sqrt R}}  \int_{|\bar x|}^{R} \tau^{-3} y^{1-2s_i} |\nabla v_i|^2 d\tau d \bar x + \frac{1}{R^2}  \int_{B^+_{R}\setminus B^+_{\sqrt R}} y^{1-2s_i} |\nabla v_i|^2 d \bar x
\\&\le & \nonumber 2 \int_{\sqrt R}^{R}  \tau^{-3}  \int_{B^+_{\tau}} y^{1-2s_i} |\nabla v_i|^2  d \bar x d\tau + \frac{1}{ R^2}  \int_{B^+_{R}} y^{1-2s_i} |\nabla v_i|^2 d \bar x.
\end{eqnarray}
Combining this and (\ref{mm}),  we conclude 
\begin{equation}\label{logR}
\int_{B^+_{R}\setminus B^+_{\sqrt R}}  \frac{y^{1-2s_i} |\nabla v_i|^2}{|\bar x|^2} d \bar x
\le C \log R,
\end{equation}
where the positive constant $C$  is independent from $R$.  % Substituting the given test function $\eta_k$ in the Poincar\'{e} geometric inequality, Theorem (\ref{lempoin}), we get 
   The above decay estimates (\ref{inelogrr}) and (\ref{logR}) imply that $|\nabla v_i|^2 \mathcal{A}_i^2 + | \nabla_{T_i} |\nabla_{ x}  v_i| |^2$ and $|\partial_{j} H_i( v)|  |\nabla_{ x}  v_i|  |\nabla_{ x}  v_j|    - \partial_{j}H_i(v)  \nabla_{ x}  v_i \cdot   \nabla_{ x}  v_j$ vanish on $\mathbb R^{3}_+$ and $\partial \mathbb R^{3}_+$, respectively. Therefore,  each $v_i$ is a one-dimensional function  and  vectors $\nabla_{ x}  v_i (x,0)$ and $   \nabla_{ x}  v_j (x,0)$ are parallel  and they are in  the same direction when $\partial_{j}H_i(v)$ is positive and in opposite directions when $\partial_{j}H_i(v)$ is negative.   This proves the desired result in two dimensions. 

We now suppose that $n=3$ and $v$ is a $H$-monotone solution of (\ref{emainh}). Note that $H$-monotonicity implies stability. The fact that $v=(v_i)_{i=1}^m$ is a bounded stable solution of (\ref{emainh}) in $\mathbb{R}^4_+$ implies that the function $\breve {v}=(\breve {v_i})_{i=1}^m$ where $\breve {v_i}(x_1, x_2,y):=\lim_{x_3\to \infty}  v_i(x_1, x_2, x_3,y)$ is also a bounded stable solution for (\ref{emainh}) in $\mathbb{R}^3_+$.  From our previous arguments regarding  $\mathbb{R}^3_+$,  we conclude reduction of dimension for each $\breve {v_i}$ that is $\breve {v_i}(x,y)=\breve {v}^*_i(\Gamma_i\cdot x,y)$ for $(x,y)\in\mathbb R^{3}_+$ and for some $\Gamma_i\in S^{1}$. From this and  Lemma \ref{asymp},  we conclude that the energy of $\breve { v}$ in $C_R\subset \mathbb{R}^3_+ $ is bounded by
 \begin{equation}
\label{boundEt} 
  E_R( \breve{ v}) \le C R^{2}, 
 \end{equation}
 when   $ E_R( \breve{ v})=\sum_{i=1}^{m}\frac{1}{2d_{s_i}} \int_{C_R}  y^{1-2s_i} |\nabla \breve {v_i}|^2 d \bar x  - \int_{B_R} \left( H( \breve{ v} ) -c_{\breve{ v}} \right)d x$ for $c_{\breve{ v}}:=\sup  H({\breve{ v}} )$.   Set $\sigma_i:=\frac{\psi_i}{\phi_i}$ when  $\phi_i := \partial _n v_i$ and $\psi_i:=\nabla v_i\cdot\eta$ for any fixed $\eta=(\eta',0)\in \mathbb{R}^{n-1}\times\{0\}$. Note that $\sigma=(\sigma_i)_{i=1}^m$ satisfies (\ref{sigma1}) and (\ref{sigma2}). To  apply Theorem \ref{thsymv},  we only need to prove the following energy estimate
\begin{equation}\label{eboundzz}
\sum_{i=1}^{m}\frac{1}{2d_{s_i}} \int_{C_R}  y^{a_i} |\nabla { v_i}|^2 d \bar x  \le  C R^{2}  \ \chi_{\{s_*>1/2\}} + C  R^{2}\log R\  \chi_{\{s_*=1/2\}} .
 \end{equation}
  Define the sequence of functions $ v^t=( v_i^t)_{i=1}^m$ when $  v_i^t(x,y):= v_i(x',x_n+t,y)$ for $t\in\mathbb{R}$ and $(x,y)=(x',x_n,y)\in\mathbb R^{4}_+$. Note that $ v^t$ is  a bounded solution of (\ref{emainh}), i.e.,
\begin{eqnarray}\label{emaint}
 \left\{ \begin{array}{lcl}
\hfill \div(y^{1-2s_i} \nabla  v^t_i)&=& 0   \ \ \text{in}\ \ \mathbb{R}_+^{n+1},\\   
\hfill -\lim_{y\to0}y^{1-2s_i} \partial_{y}  v^t_i&=& d_{s_i}  H_i( v^t)   \ \ \text{in}\ \ \partial\mathbb {R}_+^{n+1}.
\end{array}\right.
  \end{eqnarray}
Straightforward calculations show that  $ v_i^t$ converges to $ v_i$ in $C^1_{loc}(\mathbb{R}^n)$ for all $i=1,\cdots,m$ and 
  \begin{equation} \lim_{t\to\infty} E_R( v^t)= E_R( v).
  \end{equation}
On the other hand, for every $T>0$ we have 
\begin{equation}\label{bbbb}
 E_R( v)= E_R( v^T)-\int_0^T \partial_t E_R( v^t) dt.
 \end{equation}
Note that applying similar arguments,  as the ones given in the proof of Theorem \ref{energylayer}, one can get a lower bound for $\partial_t E_R( v^t)$ of the form (\ref{tEt}).   From this and  (\ref{bbbb}), we have 
%Substituting an upper bound on $\partial_t E_R(\breve v^t)$ in (\ref{bbbb}) by  following similar ideas and techniques as Theorem \ref{energylayer} we get the following upper bound on the energy that is  
 \begin{equation}
E_R( v) \le  E_R( v^T) + M \int_{\partial B_R}  \int_0^R  \sum_{i=1}^{m} \frac{y^{1-2s_i} }{1+y} dy d\mathcal H^{n-1} + M  \int_{B_R\times\{y=R\}}  \sum_{i=1}^{m} y^{-2s_i} dx.
   \end{equation}
Now, taking a limit when  $T\to\infty$ and using (\ref{boundEt}) we get 
    \begin{eqnarray}\label{ervv}
E_R( v) &\le&  E_R(\breve{ v}) + M \int_{\partial B_R}  \int_0^R  \sum_{i=1}^{m} \frac{y^{1-2s_i} }{1+y} dy d\mathcal H^{n-1} + M  \sum_{i=1}^{m} R^{n-2s_i} \\
&\le & \nonumber M R^2 + M \int_{\partial B_R}  \int_0^R  \sum_{i=1}^{m} \frac{y^{1-2s_i} }{1+y} dy d\mathcal H^{n-1} + M  \sum_{i=1}^{m} R^{n-2s_i} 
\\&\le&\nonumber  M R^2 +M  \sum_{i=1}^{m}  R^{n-2s_i} \chi_{\{0<s_i<1/2\}} + M  \sum_{i=1}^{m} R^{n-1} \chi_{\{1/2<s_i<1\}} \\&& \nonumber +  M R^{n-1}\log R \ \chi_{\{s_i=1/2\}}   +   M  \sum_{i=1}^{m} R^{n-2s_i}  . 
   \end{eqnarray}  
Here, we have used estimates provided in Lemma \ref{asymp}. Note that when $1/2 < s_*<1$ and $n=3$,  for any $1\le i\le m$,  we have 
\begin{equation}
n-2s_i\le n-2s_* < n-1.\end{equation} 
Applying (\ref{ervv}) when $n=3$ and $1/2 < s_*<1$, we conclude 
\begin{equation}
E_R( v)  \le  M R^2 +M R^{n-1}  \le C R^2,
\end{equation}
when $C$ is a positive constant that is independent from $R$.  Similarly, for the case of $s_*=1/2$ and $n=3$,  we have 
\begin{equation}
E_R( v)  \le  M R^2 +M R^{n-1} \log R \le C R^2 \log R.\end{equation}
This completes the the proof. 

               \hfill $ \Box$
  
Methods and  Ideas  applied in the above proof are strongly motivated by the ones given by Ambrosio and Cabr\'{e} in \cite{ac},  by Alberti, Ambrosio and Cabr\'{e} in  \cite{aac}, by Farina, Scuinzi and Valdinoci in \cite{fsv}, by Ghoussoub and the author in \cite{fg} and Sire and the author in \cite{fs}.     We now provide another consequence of Theorem   \ref{liouville}. The following theorem clarifies the behaviour of derivatives of each $v_i$ in various directions when  $v=(v_i)_{i=1}^m$ is a bounded stable solution of symmetric system (\ref{emainh}) in lower dimensions.       
               
\begin{thm}\label{stablev}
Suppose that $ v=(v_i)_{i=1}^m$ is a bounded stable solution of symmetric system (\ref{emainh}) for $n \le 2$ when  $s_*$ belong to $[\frac{1}{2},1)$ and for $n \le 1+2s_*$  when  $s_*$ belongs to $(0,\frac{1}{2})$. Then,  either each $\partial_x v_i(x,y)$ vanishes in $\overline{\mathbb{R}^{n+1}_+}$ or it does not change sign in $\overline{\mathbb{R}^{n+1}_+}$ for every $1\le i \le m$. 
\end{thm}
\noindent\textbf{Proof:} Let $v$ be a bounded stable solution of  (\ref{emainh}). Then, there exits a sequence of functions $\mathbf\phi=(\phi_i)_i$ such that  each $\phi_i$ does not change sign and it satisfies  
  \begin{eqnarray}\label{phii}
 \left\{ \begin{array}{lcl}
\hfill \div(y^{a_i} \nabla \phi_i)&=& 0   \ \ \text{in}\ \ \mathbb{R}_+^{n+1},\\   
\hfill -\lim_{y\to0}y^{a_i} \partial_{y} \phi_i&=& d_{s_i} \sum_{j=1}^m \partial_{j}H_i( v) \phi_j   \ \ \text{in}\ \ \partial\mathbb {R}_+^{n+1}. 
\end{array}\right.
  \end{eqnarray}
For  each index $1\le i \le m$,  define the quotient  $\sigma_i:=\frac{\partial_x v_i}{\phi_i}$ that implies $(\sigma_i\phi_i)^2=(\partial_x v_i)^2$. From Lemma \ref{deltavx},  we have 
\begin{equation}
   \int_{C_R}  y^{a_i} (\sigma_i \phi_i)^2 d\bar x \le    \int_{C_R}  y^{a_i} |\nabla v_i|^2 d\bar x \le C R^n \int_{0}^R \frac{y^{a_i}}{1+y} dy ,
   \end{equation}
 where $R>1$. Therefore, straightforward calculations show that 
   \begin{eqnarray}\label{phii}
 \int_{C_R}  y^{a_i} (\sigma_i \phi_i)^2 d\bar x \le   C\left\{ \begin{array}{lcl}
R^{n+1-2s_i}  && \ \text{for}\ \  \ 0<s_i<\frac 1 2,\\
R^n \ln R && \  \text{for}\ \ s_i=\frac 1 2,\\   
R^n  && \ \text{for}\ \ \frac 1 2 <s_i<1. 
\end{array}\right.
  \end{eqnarray}
 Note that for $\sigma_i$ and $\phi_i$ equations (\ref{sigma1}) and (\ref{sigma2}) hold. For  $n \le 2$ when all $s_i$ belong to $[\frac{1}{2},1)$ and for $n \le 1+2s_*$  when at least one of $s_i$ belongs to $(0,\frac{1}{2})$,  estimate (\ref{liouassum})  holds for an appropriate $F\in\mathcal F$.    Set $h_{i,j}:=\partial_{j}H_i( v) \phi_i \phi_j $ and  $f$ to be the identity function in Theorem \ref{liouville}. Note that for symmetric systems,  we have $h_{i,j}=h_{j,i}$.  Theorem \ref{liouville} implies that  each $\sigma_i$ is constant.  Therefore,  there exists a sequence of constants $ C=(C_i)_i$ such that $\partial_x v_i (x,y)=C_i \phi_i(x,y)$ for $(x,y)\in\overline{\mathbb{R}^{n+1}_+}$. Since  each $\phi_i$ does not change sign, the proof is completed.  

 \hfill $ \Box$

The De Giorgi's conjecture provides a  reduction of dimensions, to one-dimension, for bounded monotone solutions   of the Lane-Emden equation when $n\le 8$.    The latter theorem provides a counterpart of the conjecture to multi-component fractional symmetric systems with a general nonlinearity.    In what follows, we assume certain extra assumptions on the sign of the nonlinearity $H$ and we establish a Liouville theorem for bounded stable solutions of (\ref{emainh}) in lower dimensions applying Theorem \ref{stablev}, Theorem \ref{liouville} and Theorem \ref{hamiltonthm}.

\begin{thm}\label{lioupositive}
Suppose that $ v=(v_i)_{i=1}^m$ is a bounded stable solution of (\ref{emainh}) when either $H_i(v) \ge0$ for all $1\le i\le m$ or $\sum_{i=1}^m v_i H_i(v)\le 0$. Then, each $v_i$ must be constant provided $n\le 2(1+s_*)$.   
\end{thm}

               \noindent\textbf{Proof:} %This is a direct consequence of Lemma \ref{Hiv}.  
 Since $v$ is a stable solution, there exists a sequence $\phi=(\phi_i)_{i=1}^m$ satisfying (\ref{stabilityh}). On the other hand, applying Lemma  \ref{Hiv}, for each $i$,  we have 
\begin{equation}\label{Hpositivel}
\int_{C_R} y^{1-2s_i} |\nabla v_i|^2 d\bar x \le C R^{n-2s_i}  .
\end{equation}
The fact that $n\le 2+2s_*$,  implies that $n-2s_i \le n-2s_* \le 2$.  Note that    $\sigma_i=\frac{\nabla v_i\cdot \eta}{\phi}$ satisfies conditions of Theorem \ref{liouville} for $F(r)=1$, $h_{i,j}=\partial_j H_i(v) \phi_i \phi_j$ and $f$ to be the identity function. Therefore, each  $\sigma_i$ must be constant for an arbitrary direction $\eta$.  This implies that there exist a constant $ \Gamma_i\in{S}^{n-1}$ and $v^*_i: \mathbb{R^+}\times  \mathbb{R}^+\to  \mathbb{R}$ such that  $v_i(x,y)=v^*_i(  \Gamma_i\cdot x,y)$ for $(x,y)\in\mathbb R^{n+1}_+$. In other words, $v=(v_i)_{i=1}^m$ is a bounded stable solution of (\ref{emainh}) for $n=1$. Applying   Theorem \ref{stablev},    we conclude  that $\partial_x v_i$ does not change sign in $\mathbb R^2_+$. Here,  we have used the fact that when $\partial_x v_i$ vanishes,  boundedness implies that each $v_i$  must be constant.  Therefore, $v_i$ has to be strictly monotone in $x$ which together with boundedness of $v_i$  proves the existence of $\lim_{x\to\pm\infty}v_i(x,0)$. Let $\lim_{x\to-\infty}v_i(x,0)=l_i$ and $\lim_{x\to\infty}v_i(x,0)=L_i$,  where $l_i$ and $L_i$ are constants. Since $v_i$ is strictly monotone, we conclude $l_i < L_i $.  We now apply the Hamiltonian identity provided in Theorem \ref{hamiltonthm} when $x\to\pm\infty$  to get $\tilde H(l)=\tilde H(L)$ where $l=(l_i)_{i=1}^m$ and $L=(L_i)_{i=1}^m$.  Note that this is in contradiction with the following 
\begin{equation}
0=\tilde H(L)-\tilde H(l)=\sum_{i=1}^m (L_i-l_i)H_i(t(L-l)+l),
\end{equation}
for some $t\in (0,1)$. This completes the proof.

               \hfill $ \Box$

 Mathematical techniques and ideas that we applied in the above proof are strongly motivated by the ones given  in \cite{df} by Dupaigne and Farina. In \cite{df}, authors proved that any  bounded stable solution of (\ref{sceqH}), that is when $m=1$ and $s=1$ in (\ref{mainh}),  is constant provided $ n \le 4$ and $0 \le H \in C^1(\mathbb R)$ is a general nonlinearity. %We refer interested readers to \cite{mf} for more information .  
%for the case of scalar equation, that is when $m=1$,  and for the local case, that is when $s=1$,  Dupaigne and Farina in \cite{df} proved that any classical bounded stable solution of (\ref{emain}) is constant provided $ n \le 4$ and $0 \le H \in C^1(\mathbb R)$ is a general nonlinearity. 
Note that for particular nonlinearities the critical dimension is much higher than four dimensions. We also refer interested readers to  \cite{cc} by Cabr\'{e} and Capella, to  \cite{sv} by Villegas and to \cite{mf} by the author for the case of radial stable solutions where the optimal dimension is $n=10$ for a general nonlinearity $H\in C^1(\mathbb R)$.   So,  we expect that Theorem \ref{lioupositive} could be improved.  For specific nonlinearities $H(u)=e^u$, $H(u)=u^p$ where $p>1$ and $H(u)=-u^{-p}$ where $p>0$ the equation is called Gelfand, Lane-Emden and Lane-Emden with negative exponent equations, and Lioville theorems are given for the following optimal dimensions, respectively, 
\begin{itemize}
\item  $1\le n <10$ by Farina in \cite{f2}, 
\item  $1\le n < 2+ \frac{4}{p-1}(p+\sqrt{p(p-1)})$  by Farina in \cite{f1}, 
\item  $1\le n < 2+ \frac{4}{p+1}(p+\sqrt{p(p+1)})$ by Esposito, Ghoussoub and Guo in \cite{egg}. 
\end{itemize}
%Note that for particular nonlinearities the critical dimension is much higher than four dimensions. So one might expect to improve the dimension for the case of general nonlinearity.  
%Let us mention that for radial solutions an optimal Liouville theorem is given  by Cabr\'{e} and Capella in  \cite{cc} and by Villegas in \cite{sv} stating that any bounded radial stable solution of $ -\Delta u= H(u)$ has to be constant provided $1\le n < 10$ for a general nonlinearity $H\in C^1(\mathbb R)$. 
%For the rest of this section we shall focus on proofs of main results.  %We start with the proof of the Hamiltonian identity. 
In the next section, we study a two-component nonlinear Schr\"{o}dinger  system,  that is a particular case of  (\ref{mainh}),   and we prove various Liouville theorems.

%\begin{open} Suppose that $v=(v_i)_{i=1}^m$ is a  bounded $H-$monotone solution of  orientable symmetric (\ref{emainh}) with a general nonlinearity $H=(H_i)_{i=1}^m$ in three dimensions and $0<s_i<1$ where $1/2\le s_* <1$ does not hold anymore.   Does the conclusion of Theorem \ref{thsymv} hold?\end{open}

%%%%%%%%%%%%%%%%%%%%%%%%%%
%%%%%%%%%%%%%%%%%%%%%%%%
\section{Nonlinear Schr\"{o}dinger  System; Particular Nonlinearity}\label{sch}
Consider the following two-component system
\begin{eqnarray}\label{mainsch}
 \left\{ \begin{array}{lcl}
\hfill (- \Delta)^{s_1}  u_1 &=& \mu_1 u_1^3+\beta u_2^2 u_1   \ \ \text{in}\ \ \mathbb{R}^n,\\   
\hfill (- \Delta)^{s_2} u_2 &=& \mu_2 u_2^3+\beta u_1^2 u_2     \ \ \text{in}\ \ \mathbb{R}^n,
\end{array}\right.
  \end{eqnarray}
when $0<s_1,s_2<1$ and $\mu_1,\mu_2,\beta$ are parameters.  This is a special case of system (\ref{mainh}) when $m=2$ and $H_1(u_1,u_2)=\mu_1 u_1^3+\beta u_2^2 u_1$ and $H_2(u_1,u_2)=\mu_2 u_2^3+\beta u_1^2 u_2 $.   The above  system arises in Bose-Einstein condensations  and it is well-studied in the literature. We refer interested readers to \cite{ww, fs,  tvz, nttv, wwe2,bdw,wwe, lw, dww} and references therein for more information. % and for a particular case of $\mu_1=\mu_2=0$ and $\beta=-1$ for both local and nonlocal cases. 
The extension function pair $(v_1,v_2)$ given in (\ref{emainh}) satisfies 
 \begin{eqnarray}\label{emain}
 \left\{ \begin{array}{lcl}
\hfill \div(y^{a_i} \nabla v_i)&=& 0   \ \ \text{in}\ \ \mathbb{R}_+^{n+1},\\  
\hfill -d_{s_i} \lim_{y\to0}y^{a_i} \partial_{y} v_i&=&  H_i(v_1,v_2)  \ \ \text{in}\ \ \partial\mathbb{R}_+^{n+1},\\
\end{array}\right.
  \end{eqnarray}
when $a_i=1-2s_i$ and $d_{s_i} = \frac{2^{2s_i-1}\Gamma(s_i)}{\Gamma(1-s_i)}$ for $1\le i\le 2$. % Here $u_e,v_e$ are the extensions of the functions $u,v$, respectively. % Main results of the present paper deal with  solutions of \eqref{emain}. The direct consequences of these results, in the light of \cite{cafs}, leads us to similar results for the original system \eqref{main}. 
Note that when $\beta=0$, system (\ref{mainsch}) becomes decoupled and each equation in (\ref{mainsch}) is of the from 
\begin{equation}\label{lanem}
 (- \Delta)^{\alpha}  w =  w^3   \ \ \text{in}\ \ \mathbb{R}^n. 
\end{equation} 
when $\mu_1,\mu_2>0$ for either $w=\sqrt{\mu_1} u_1$ and $\alpha=s_1$ or $w=\sqrt{\mu_2} u_2$ and $\alpha=s_2$. The above equation (\ref{lanem}) is  known as the fractional Lane-Emden equation and nonnegative solutions of this equation are classified completely in the literature.   It is known that whenever 
\begin{equation}\label{n4alpha}
n<4\alpha,
\end{equation}  the only nonnegative solution for (\ref{lanem}) is the trivial solution and $n=4\alpha $ is the critical dimension, see \cite{li, gs, clo}.  We refer interested readers to \cite{ddw,fw,ddww} for the classification of stable solutions of (\ref{lanem}) when $0<\alpha\le 2$.  In this article, we are interested in the case of $\beta\neq 0$. % and the more challenging case that is when $\beta<0$.  
   The following Liouville theorem addresses stable solutions of (\ref{emain}) and is a direct consequence of Theorem \ref{lioupositive}.  The proof is straightforward and we omit it here. 

%The next result is a Liouville theorem for stable solutions. Note that positivity of solutions is not assumed for Part (i). 

\begin{thm}\label{thmschstable}
Let $v=(v_1,v_2)$ be a bounded stable solution of (\ref{emain}) when $n\le 2+2 \min\{s_1,s_2\}$. Assume that   either $\mu_1,\mu_2\le 0$  and $|\beta| \le \sqrt{\mu_1\mu_2}$
or $\mu_1,\mu_2,\beta \ge 0$ and solutions $(v_1,v_2)$ are nonnegative.  Then,  each $v_i \equiv C_i$  where $C_i$ is constant.  
\end{thm}

%\noindent\textbf{Proof:} We omit the proof. %since it is a direct consequence of Theorem \ref{lioupositive}.  \hfill $ \Box$

We now provide a Liouville theorem for solutions of (\ref{mainsch}) in the absence of stability.

\begin{thm}\label{thmsch}
Suppose that  $s_1=s_2=\alpha$ and $\mu_1,\mu_2$ are nonnegative. Assume that  $u=(u_1,u_2)$ is a nonnegative solution of (\ref{mainsch}) when $n\le 3\alpha$ and $\beta >-\sqrt{\mu_1\mu_2}$, then $(u_1,u_2)=(0,0)$. 
\end{thm}

\noindent\textbf{Proof:}  Let $u=(u_1,u_2)$ be a nonnegative solution of  (\ref{mainsch}) when $s_1=s_2=\alpha$. Without loss of generality assume that $u_1>0$ and $\mu_1>0$.  It is straightforward to show that equivalently $(u_1,u_2)$ is a solution for  the integral system  
\begin{eqnarray}\label{emainin}
 \left\{ \begin{array}{lcl}
\hfill  u_1(x) &=& \int_{\mathbb R^n} \frac{\mu_1 u_1^3(y)+\beta u_2^2(y) u_1(y)}{|x-y|^{n-2\alpha}} dy   \ \ \text{for} \ x\in\mathbb R^n,\\   
\hfill u_2 (x) &=& \int_{\mathbb R^n} \frac{\mu_2 u_2^3(y)+\beta u_1^2(y) u_2(y)}{|x-y|^{n-2\alpha}} dy       \ \ \text{for} \ x\in\mathbb R^n
\end{array}\right.
  \end{eqnarray}
Assume that parameter $\beta$ is nonnegative.  Since each $u_i$ is nonnegative, system (\ref{emainin})  gives  
\begin{eqnarray}\label{emaininb}
 \left\{ \begin{array}{lcl}
\hfill  u_1(x) &\ge&\mu_1  \int_{\mathbb R^n} \frac{ u_1^3(y)}{|x-y|^{n-2\alpha}} dy    \ \ \text{for} \ x\in\mathbb R^n,   \\ 
\hfill  u_2(x) &\ge &\mu_2  \int_{\mathbb R^n} \frac{ u_2^3(y)}{|x-y|^{n-2\alpha}} dy       \ \ \text{for} \ x\in\mathbb R^n.
\end{array}\right.
  \end{eqnarray}
Applying Liouville theorems given in \cite{cam1}, see also \cite{cam2,clo,li}, to each of the above integral inequalities we conclude that  $(u_1,u_2)=(0,0)$ when $n\le 3 \alpha$. We now let $\beta $   be negative that is $-\sqrt{\mu_1 \mu_2} < \beta <0$. Multiplying the first equation  in (\ref{emainin}) with $\sqrt[4]{{\mu_2}/{\mu_1}} $ and adding both equations in (\ref{emainin}) we  get % for $x\in\mathbb R^n$
\begin{eqnarray}\label{4mu1mu2}
 && \sqrt[4]{{\mu_2}/{\mu_1}} u_1(x)+ u_2(x) =
 \\&& \nonumber \int_{\mathbb R^n} \frac{1}{|x-y|^{n-2\alpha}}   \left[\sqrt[4]{{\mu_2}/{\mu_1}} \left(  \mu_1 u_1^3(y)  +\beta u_1(y) u_2^2(y)\right)   +\mu_2 u_2^3(y)  +\beta u_2(y) u_1^2(y)  \right] dy, 
  \end{eqnarray}
 for $x\in\mathbb R^n$. To simplify the right-hand side of above equality, we claim that there exists a positive constant $M$ such that %for all nonnegative $u_i$, 
\begin{equation}\label{M}
\sqrt[4]{{\mu_2}/{\mu_1}} \left(  \mu_1 u_1^3  +\beta u_1 u_2^2\right)   +\mu_2 u_2^3  +\beta u_2 u_1^2  \ge M \left( \sqrt[4]{{\mu_2}/{\mu_1}} u_1+u_2\right)^3. 
\end{equation}
In order to prove this claim,  let $\tau:={u_2}/{u_1}$ and  define continuous function $F:\mathbb R^+\to \mathbb R$ as 
\begin{equation}
F(\tau)=\frac{ \sqrt[4]{{\mu_2}/{\mu_1}} \left(\mu_1+\beta \tau^2\right) + \mu_2 \tau^3 + \beta \tau}{(\sqrt[4]{ {\mu_2}/{\mu_1}}+\tau)^3 }.
  \end{equation}
 Note that $F(0)=\mu_1>0$ and $\lim_{\tau\to \infty} F(\tau)=\mu_2 \ge 0$. In addition, for $\beta>-\sqrt{\mu_1\mu_2}$,  we have
 \begin{eqnarray}
\left(  \sqrt[4]{{\mu_2}/{\mu_1}}+\tau \right)^3 F(\tau) &>&  \sqrt[4]{{\mu_2}/{\mu_1}} \left( \mu_1 - \sqrt{\mu_1\mu_2} \tau^2  \right) +\tau (\mu_2 \tau^2 - \sqrt{\mu_1\mu_2})
\\&=&  \left( \sqrt[4]{\mu_2} \tau-   \sqrt[4]{\mu_1} \right)^2  \left( \sqrt{\mu_2} \tau+  \sqrt[4]{\mu_1 \mu_2} \right) \ge  0 , 
   \end{eqnarray}
 that is $F(\tau) >0$ when $\tau>0$.  This proves the claim.  From  (\ref{M}) and (\ref{4mu1mu2}),  we get 
\begin{equation}\label{u+v}
\sqrt[4]{{\mu_2}/{\mu_1}} u_1(x)+ u_2(x) \ge M  \int_{\mathbb R^n} \frac{1}{|x-y|^{n-2\alpha}}  \left( \sqrt[4]{{\mu_2}/{\mu_1}} u_1(y)+u_2(y)  \right)^3 dy \ \ \text{for} \ \ x\in\mathbb R^n. 
  \end{equation}
Let  $z(x):= \sqrt{M} \left(\sqrt[4]{{\mu_2}/{\mu_1}} u_1(x)+ u_2(x)\right) >0$ in (\ref{u+v}) to obtain 
\begin{equation}
 z(x) \ge  \int_{\mathbb R^n} \frac{z^3(y)}{|x-y|^{n-2\alpha}} dy \ \ \text{for} \ x\in\mathbb R^n. 
  \end{equation}
%Scaling $\bar z := \sqrt{M} z$ we obtain  \begin{equation} z(x) \ge  \int_{\mathbb R^n} \frac{ z^3(y)}{|x-y|^{n-2\alpha}} dy \ \ \text{for} \ x\in\mathbb R^n.  \end{equation}
We now apply Liouville theorems given in \cite{cam1,cam2} for the above integral inequality  to conclude  that $z\equiv0$ when $n\le 2\alpha $ and  for $n>2\alpha$ when $3 \le \frac{n}{n-2\alpha}$. This implies that whenever $n\le 3\alpha$,  we have $z\equiv0$. This completes the proof. 

\hfill $ \Box$

According to (\ref{n4alpha}), one may expect that the critical dimension for the above theorem is $n=4 \alpha$ and this remains an open problem. Note that methods and ideas applied in the above proof are strongly motivated by the ones given in \cite{dww} by Dancer, Wei and Weth and in \cite{lw} by Lin and Wei. %We refer interested readers to  \cite{bdw,wwe,wwe2} for more information. *****
We now provide  a monotonicity formula for solutions of the two-component Schr\"{o}dinger system (\ref{emain}). This is a direct consequence of Theorem \ref{thmmonoI} and we omit the proof.  %for not necessarily radial solutions for some specific parameters $\mu_1,\mu_2,\beta$. 
Note that Frank and Lenzmann in \cite{fl} and with Silvestre in \cite{fls} used similar monotonicity formulae  to study uniqueness of solutions for the fractional Schr\"{o}dinger operator.

\begin{thm} 
Let $v=(v_1,v_2)$ be a bounded solution of \eqref{emain}. Suppose that $\mu_1,\mu_2\le 0$ and $|\beta | \le \sqrt{\mu_1 \mu_2}$.  Then 
\begin{equation}
I(R) =   R^{-n+2s^*} \left[ \sum_{i=1}^2 d_{s_i}   \int_{B_R^+} y^{a_i}  |\nabla v_i|^2 d\bar x  
 -  \int_{B_R\times\{y=0\}} 
\left( \frac{\mu_1}{2} v_1^4 + \frac{\mu_2}{2} v_2^4 + \beta v_1^2 v_2^2\right) d x \right] , 
\end{equation}
 is a nondecreasing function of $R \ge 1$. %for parameters $\mu_1,\mu_2,\beta $ such that $\mu_1 u_e^4 +  \mu_2 v_e^4 + 2 \beta u_e^2 v_e^2 \le 0$.    
\end{thm}  
%\begin{rem} For example, when $\mu_1,\mu_2\le 0$ and $\beta\le \sqrt{\mu_1\mu_2}$ then we have this monotonicity formula.  \end{rem}

Lastly,  we provide a monotonicity formula for radial solutions of \eqref{emain}  for all parameters $\mu_1,\mu_2,\beta\in\mathbb R$. 
%Note that understanding of qualitative behaviour of solutions of system (\ref{main})  is much more challenging, and then  more interesting,  when $\mu_1,\mu_2$ are nonnegative  and $\beta $ is negative.  In addition to Liouville theorems and monotonicity formulas provided in this paper,   we refer interested readers to \cite{ww,fs, tvz,nttv} and references therein for a particular case of $\mu_1=\mu_2=0$ and $\beta=-1$ for both local and nonlocal cases where system (\ref{main}) arises in Bose-Einstein condensation. 
%******

\begin{thm}\label{hamilton}
 Let $v=(v_i)_{i=1}^m$ be a bounded radial solution of \eqref{emain}, i.e. $v_i(x,y)=v_i(|x|,y)$. Then, the following function is nondecreasing in $r>0$,
  \begin{equation}\label{radialmono}
  J(r):= \sum_{i=1}^2 d_{s_i} \int_0^\infty   y^{1-2s_i} \left[  (\partial_r v_i(r,y))^2 - (\partial_y v_i(r,y))^2\right] dy  - \tilde H(v_1(r,0),v_2(r,0)) , 
  \end{equation}
  where $ \tilde H(v) =  \frac{\mu_1}{2} v_1^4 + \frac{\mu_2}{2} v_2^4 + \beta v_1^2v_2^2.$
 More precisely, 
  \begin{equation}\label{radialmonoprime}
  J'(r)=-\frac{n-1}{r} \sum_{i=1}^2 d_{s_i}  \int_0^\infty  y^{1-2s_i} (\partial_r v_i)^2 dy .
    \end{equation}
\end{thm}

\noindent\textbf{Proof:} Suppose that $v=(v_i)$ is radially symmetric in $ x$. For  $r=| x|$ we have
 \begin{eqnarray}\label{emainrad}
 \left\{ \begin{array}{lcl}
\hfill \partial_{rr} v_i +\frac{n-1}{r} \partial_r v_i + \partial_{yy} v_i +\frac{a_i}{y} \partial_y v_i&=& 0   \ \ \text{in}\ \  (0,\infty)\times(0,\infty),\\  
\hfill -d_{s_{i}} \lim_{y\to0}y^{a_i} \partial_{y} v_i &=&  H_i(v)   \ \ \text{in}\ \ (0,\infty)\times(y=0),\\
\end{array}\right.
  \end{eqnarray}
where $a_i=1-2s_i$.  Define the following function of $r$, 
  \begin{eqnarray}
  p(r) &:=& \sum_{i=1}^2 \frac{d_{s_i}}{2} \int_0^\infty  y^{1-2s_i} \left[  (\partial_r v_i)^2 - (\partial_y v_i)^2\right]dy 
    \end{eqnarray}
 Taking derivative of $p$ with respect to $r$ and using (\ref{emainrad}) to  substitute values of $ \partial_{rr} v_i$,  we conclude 
 \begin{eqnarray}\label{wr}
 p'(r) &=&   -\frac{n-1}{r} \sum_{i=1}^2  d_{s_i} \int_0^\infty   y^{a_i} (\partial_r v_i)^2 dy  -   \sum_{i=1}^2  d_{s_i}   \int_0^\infty   y^{a_i}  ( \partial_{yy} v_i)(  \partial_{r} v_i )  dy
 \\&& \nonumber -   \sum_{i=1}^2  d_{s_i}  a_i   \int_0^\infty   y^{-2s_i} (\partial_{r} v_i) (\partial_{y} v_i)    dy 
-   \sum_{i=1}^2  d_{s_i}   \int_0^\infty   y^{a_i} (\partial_{y} v_i) (\partial_{ry} v_i)  dy  .
  \end{eqnarray}
   Applying integration by parts and using the boundary term in (\ref{emainh}),   we have   
  \begin{eqnarray}\label{intpart}
-\sum_{i=1}^2  d_{s_i}  \int_0^\infty  \left[ y^{a_i}  \partial_{yy} v_i \partial_{r} v_i   + a_i    y^{-2s_i} \partial_{r} v_i \partial_{y} v_i +      y^{a_i}  \partial_{y} v_i  \partial_{ry} v_i \right]dy
 & =& \lim_{y\to 0}\sum_{i=1}^2  d_{s_i}   y^{a_i}  \partial_{r} v_i \partial_{y} v_i 
  \\&=&   \nonumber \partial_{r}\left(\tilde H(v) \right) . 
   \end{eqnarray}
   Combining (\ref{intpart}) and (\ref{wr}) and setting $J(r):=p(r)-\tilde H(v(r,0))$ completes the proof.

               \hfill $ \Box$

We end this section with this point that in the absence of monotonicity and stability assumptions,  the qualitative behaviour of solutions of elliptic and Hamiltonian systems with a general nonlinearity are studies extensively in the literature. We  refer interested readers to \cite{dff, dfd, gui} and references therein.


\begin{thebibliography}{99}   

\bibitem {aac} G. Alberti, L. Ambrosio, and X. Cabr\'{e}, \emph{On a long-standing conjecture of E. De Giorgi: symmetry in 3D for general nonlinearities and a local minimality property}, Acta Appl. Math. 65 (2001), 9-33. 

\bibitem{ac} L. Ambrosio and X. Cabr\'{e}, \emph{Entire solutions of semilinear elliptic equations in $\mathbb R^3$ and a conjecture of De Giorgi}, J. Amer. Math. Soc. 13 (2000), 725-739. 


\bibitem{bdw} T. Bartsch, N. Dancer, Z. Q. Wang,  \emph{A Liouville theorem, a-priori bounds, and bifurcating
branches of positive solutions for a nonlinear elliptic
system}, Calc. Var. (2010) 37:345-361. 

\bibitem{bcn} H. Berestycki, L. Caffarelli, L. Nirenberg,  \emph{Further qualitative properties for elliptic equations in unbounded domains}. Ann. Scuola Norm. Sup. Pisa Cl. Sci. 25 (1997) 69-94.

\bibitem{cc}  X. Cabr\'{e}, A. Capella, \emph{On the stability of radial solutions of semi-linear elliptic equations in all of $R^n$}, C. R. Math. Acad. Sci. Paris 338 (2004) 769-774.

\bibitem{ccinti1} X. Cabr\'{e},  E.  Cinti, \emph{Energy estimates and 1-D symmetry for nonlinear equations involving the half-Laplacian}. Discrete Contin. Dyn. Syst. 28 (2010), no. 3, 1179-1206.

\bibitem{ccinti} X. Cabr\'{e},  E.  Cinti, \emph{Sharp energy estimates for nonlinear fractional
diffusion equation}. Calc. Var. PDEs (2014) 49:233-269.

\bibitem{cs1}  X. Cabr\'{e}, Y. Sire, \emph{Nonlinear equations for fractional Laplacians I:
Regularity, maximum principles, and Hamiltonian estimates}. Annales De l'IHP,
Analyse Nonlinéaire (2014) 31:23-53.

\bibitem{cs2}  X. Cabr\'{e}, Y. Sire, \emph{Nonlinear equations for fractional Laplacians II:
Existence, uniqueness and qualitative properties of solutions}. Trans. AMS
(2015) 367:911-941.



\bibitem{cso}  X. Cabr\'{e},  J. Sol\'{a}-Morales, \emph{Layer solutions in a half-space for boundary reactions}, Comm. Pure and Appl. Math. 58 (2005), 1678-1732

\bibitem{cafs} L. Caffarelli, L. Silvestre, \emph{An extension problem related to the
fractional Laplacian}. Commun. PDE (2007) 32:1245-1260.

\bibitem{cam1} G. Caristi, L. D'Ambrosio, E. Mitidieri, \emph{Liouville theorems for some nonlinear inequalities}, Proceedings of the Steklov Institute of Mathematics, (2008), Vol. 260, pp. 90-111. 

\bibitem{cam2} G. Caristi, L. D'Ambrosio, E. Mitidieri, \emph{Representation formulae for solutions to some classes of higher order systems and related Liouville theorems},
 Milan Journal of Mathematics 76 (2008), 27-67.  
 
 \bibitem{clo} W. Chen, C. Li, and B. Ou, \emph{Classification of solutions for an integral equation}, Comm. Pure Appl. Math. 59:3 (2006), 330-343.


\bibitem{dww}  EN. Dancer, J. Wei, T. Weth, \emph{A priori bounds versus multiple existence of positive solutions for a nonlinear Schr\"{o}dinger system}, AIHP 27 (2010), pp. 953-969. 

    \bibitem{ddw} J. Davila, L. Dupaigne, J. Wei, \emph{On the fractional Lane-Emden equation}, To appear in Trans. Amer. Math. Soc.

 \bibitem{ddww} J. Davila, L. Dupaigne,  K. Wang, J. Wei, \emph{A monotonicity formula and a Liouville-type theorem for a fourth order supercritical problem}, Advances in Mathematics 258 (2014), 240-285.
 
\bibitem{dfd} D. G. De Figueiredo,  Y. H. Ding, \emph{Strongly indefinite functionals and multiple solutions of elliptic systems}, Trans. Amer. Math. Soc. 355 (2003), 2973-2989. 

\bibitem{dff} D. G. De Figueiredo,  P. L. Felmer, \emph{On superquadratic elliptic systems}, Trans. Amer. Math. Soc. 343 (1994), no. 1, 99-116.  


\bibitem{DeGiorgi} E. De Giorgi, \emph{Convergence problems for functional and operators}, Proceedings of the International Meeting on Recent Methods in Nonlinear Analysis (Rome, 1978), pp. 131-188, Pitagora, Bologna, 1979. 

\bibitem{dkw} M. del Pino, M. Kowalczyk, J. Wei,   \emph{On De Giorgi's conjecture in
dimension $N\ge 9$}. Ann. of Math. (2) (2011). 174:1485-1569. 

\bibitem{dp} S. Dipierro, A.  Pinamonti, \emph{A geometric inequality and a symmetry result for elliptic systems involving the fractional Laplacian}. J. Differential Equations 255 (2013), no. 1, 85-119.
 
 \bibitem{df} L. Dupaigne, A. Farina, \emph{Stable solutions of $-\Delta u = f(u)$ in $R^n$}, J. Eur.
Math. Soc. (JEMS) (2010) 855-882.

\bibitem{ds} L. Dupaigne, Y. Sire, \emph{A Liouville theorem for non local elliptic
equations}, Contemporary Mathematics, Vol.
528. Providence, RI: AMS, (2010) pp. 105-114. 

\bibitem{egg} P. Esposito, N. Ghoussoub,  Y. Guo, \emph{Compactness along the branch of semistable and unstable solutions for an elliptic problem with a singular nonlinearity},  Comm. Pure Appl. Math. 60 (2007), no. 12, 1731-1768. 


\bibitem{f1} A. Farina, \emph{On the classification of solutions of the Lane-Emden equation on unbounded domains of $\mathbb R^n$}, J. Math. Pures Appl. (9) 87 (2007), no. 5, 537-561.

\bibitem{f2} A. Farina, \emph{Stable solutions of $-\Delta u = e^u$ on $\mathbb R^n$}, C. R. Math. Acad. Sci. Paris 345 (2007), no. 2, 63-66. 

\bibitem{fsv} A. Farina, B. Sciunzi, E.  Valdinoci, \emph{Bernstein and de giorgi type
problems: New results via a geometric approach}. Ann. Sc. Norm. Super. Pisa Cl.
Sci. (2008)  7:741-791.

\bibitem{fn} A. Farina, N.  Soave,  \emph{Monotonicity and 1-dimensional symmetry for solutions of an elliptic system arising in Bose-Einstein condensation}. Arch. Ration. Mech. Anal. 213 (2014), no. 1, 287-326. 

\bibitem{mf} M. Fazly, \emph{Rigidity results for stable solutions of symmetric systems}, 
Proc. Amer. Math. Soc. 143 (2015), 5307-5321. 

\bibitem{fg} M. Fazly,  N. Ghoussoub, \emph{De Giorgi type results for elliptic systems}, 
Calc. Var. Partial Differential Equations 47 (2013) 809-823.

\bibitem{fs} M. Fazly, Y. Sire, \emph{Symmetry results for fractional elliptic systems and related problems}, Communications in PDEs, 40 (2015) 1070-1095. 


\bibitem{fw}M. Fazly, J. Wei, \emph{On finite Morse index solutions of higher order fractional Lane-Emden equations}, To appear in  Amer. J. Math.


 \bibitem{fl} R.L. Frank and E. Lenzmann, \emph{Uniqueness of non-linear ground states for fractional Laplacians in $\mathbb R$}, Acta Math. 210 (2013), no.2, 261-318.
 
\bibitem{fls}  R.L. Frank, E. Lenzmann and L. Silvestre, \emph{Uniqueness of radial solutions for the fractional Laplacian},   Comm. Pure Appl. Math. (2016)  69 (9). pp. 1671-1726. 


\bibitem {gg} N. Ghoussoub, C. Gui, \emph{On a conjecture of De Giorgi and some related problems}, Math. Ann. 311 (1998), no. 3, 481-491.

\bibitem {gg2} N. Ghoussoub,  C. Gui,  \emph{On De Giorgi's conjecture in dimensions 4 and 5}. Ann. Math. (2) 157(1), 313-334 (2003)

\bibitem{gs}  B. Gidas, J. Spruck, \emph{{A priori bounds for positive solutions of nonlinear elliptic equations}, Comm. Partial Differential Equations}  6 (1981) 883-901.

\bibitem{gui} C. Gui,  \emph{Hamiltonian identities for elliptic partial differential equations}.
J. Funct. Anal. (254)  (2008) 904-933.

\bibitem{k1}  L. Karp, \emph{Asymptotic behavior of solutions of elliptic equations I: Liouville-type theorems for linear and nonlinear equations on $\mathbb R^n$}, J. Analyse Math. 39 (1981) 75-102.

\bibitem{k2} L. Karp, \emph{Asymptotic behavior of solutions of elliptic equations II: Analogues of Liouville's theorem for solutions of inequalities on $\mathbb  R^n$, $n  \ge 3$}, J. Analyse Math. 39 (1981) 103-115.

\bibitem{li} Y. Y. Li, \emph{Remark on some conformally invariant integral equations: The method of moving spheres}, J. Eur. Math. Soc. (JEMS) 6:2 (2004), 153-180. 

\bibitem{lw}  T.C. Lin, J. Wei, \emph{Ground state of $N$ coupled nonlinear Schr\"{o}dinger equations in $\mathbb R^n$,  $n \le 3$}, Comm. Math. Phys. 255:3 (2005), 629-653.

\bibitem{mod} L. Modica, \emph{A gradient bound and a Liouville theorem for nonlinear Poisson equations}, Comm. Pure Appl. Math. 38 (1985), 679-684.

\bibitem{mos} L. Moschini, \emph{New Liouville theorems for linear second order degenerate elliptic equations in divergence form}, Ann. I. H. Poincar\'{e} - AN 22 (2005) 11-23.


\bibitem{nttv} B. Noris, H. Tavares, S. Terracini, G. Verzini, \emph{Uniform H\"{o}lder bounds for nonlinear Schr\"{o}dinger systems with strong competition}, Comm. Pure Appl. Math. 63 (2010), 267-302. 

\bibitem{sav} O. Savin,  \emph{Regularity of flat level sets in phase transitions}. Ann. of Math.
(2)  (2009) 169:41-78. 


\bibitem{LS} L. Silvestre, \emph{Regularity of the obstacle problem for a fractional power of the Laplace operator}, Comm. Pure Appl. Math. 60 (2007), no. 1, 67-112.

\bibitem{sv} Y. Sire, E. Valdinoci, \emph{Fractional Laplacian phase transitions and boundary reactions: A geometric inequality and a symmetry result}, J. Functional Analysis 256 (2009), 1842-1864.

\bibitem{sz}  P. Sternberg, K. Zumbrun, \emph{A Poincar\'{e} inequality with applications to
volume-constrained area-minimizing surfaces}. J. Reine Angew. Math.  (1998) 503:63-85. 

\bibitem{tvz} S. Terracini, G. Verzini, A.  Zilio, \emph{Uniform H\"{o}lder regularity with small exponent in competition-fractional diffusion systems}, Discrete Contin. Dyn. Syst. 34 (2014), no. 6, 2669-2691.  


\bibitem{ww} K. Wang, J. Wei, \emph{On the uniqueness of solutions of an nonlocal elliptic system},  Math. Annalen 365(2016), no.1-2, 105-153. 

    
\bibitem{wwe} J. Wei, T. Weth, \emph{Radial solutions and phase separation
in a system of two coupled Schr\"{o}dinger equations}, Arch. Rational Mech. Anal. 190 (2008) 83-106. 

\bibitem{wwe2} J. Wei, T. Weth, \emph{Nonradial symmetric bound states for a system of
coupled Schr\"{o}dinger equations}, Rend. Lincei Mat. Appl. 18:3 (2007), 279-293. 

\end{thebibliography}
       \end{document}